\documentclass{amsart}
\begin{document}

% ----------------------------------------------------------------
\vfuzz2pt % Don't report over-full v-boxes if over-edge is small
\hfuzz2pt % Don't report over-full h-boxes if over-edge is small
%Theorems
\newtheorem{thm}{Theorem}[section]
\newtheorem{corollary}[thm]{Corollary}
\newtheorem{lemma}[thm]{Lemma}
\newtheorem{proposition}[thm]{Proposition}

\newtheorem{defn}[thm]{Definition}

\newtheorem{remark}[thm]{Remark}
\newtheorem{example}[thm]{Example}
\newtheorem{fact}[thm]{Fact}

\def\proof{\medskip Proof.\ }
\font\lasek=lasy10 \chardef\kwadrat="32 %kwadrat
\def\kwadracik{{\lasek\kwadrat}}
\def\koniec{\hfill\lower 2pt\hbox{\kwadracik}\medskip}

\def\R{\mathbb{R}}
\def\C{\mathbb{C}}
\def\T{\mathit{T}}

\def\det{\hbox{\rm det}\, }
\def\detc{\hbox{\rm det }_{\C}}
\def\i{\hbox{\rm i}}
\def\tr{\hbox{\rm tr}\, }
\def\rk{\hbox{\rm rk}\,}
\def\vol{\hbox{\rm vol}\,}
\def\Im {\hbox{\rm Im}\, }
\def\Re{\hbox{\rm Re}\, }
\def\interior{\hbox{\rm int}\, }
\def\e{\hbox{\rm e}}
\def\pu{\partial _u}
\def\pv{\partial _v}
\def\pui{\partial _{u_i}}
\def\puj{\partial _{u_j}}
\def\puk{\partial {u_k}}
\def\div{\hbox{\rm div}\,}
\def\Ric{\hbox{\rm Ric}\,}
\def\ker{\hbox{\rm ker}\,}
\def\im{\hbox{\rm im}\, }
\def\I{\hbox{\rm I}\,}
\def\id{\hbox{\rm id}\,}
\def\exp{\hbox{{\rm exp}^{\tilde\nabla}}\.}
\def\M{\mathcal M}
\def\tRic{\Ric^{\tg}}

\def\gnu{\nu_g}
\def\stat{(g,\nabla)}
\def\statstar{(\nabla ^*,g)}
\def\gnabla{\nabla^g}
\def\nablastar{\nabla^*}
\def\sup{\hbox{\rm sup}}
\def\gRic{{\Ric^g}}
\def\cinfty{\mathcal C^\infty}
\def\Rstar{R^*}
\def\gR{R^g}
\def\Ricstar{\Ric^*}
\def\RicK{\Ric ^K}
\def\rhog{\rho ^g}
\def\kK{k^K}
\def\knabla{k^\nabla}
\def\S{\mathcal S}
\def\K{\mathfrak K}
\def\tnabla{\nabla^{\tilde g}}
\def\tg{\tilde g}
\def\tR{R^{\tg}}
\def\hX{{X^h}}
\def\vX{{X^v}}
\def\hY{{Y^h}}
\def\vY{{Y^v}}
\def\hZ{{Z^h}}
\def\vZ{{Z^v}}
\def\hW{{W^h}}
\def\vW{{W^v}}
\def\F{\mathcal F}
\def\u{\xi}
\def\tk{k^{\tg}}
\def\rhotg{\rho^{\tg}}
\def\gk{k^g}
\def\Kk{k^K}

\def\U{\mathcal U}
\def\Uv{\mathcal U}
\def\Sr{S^r\M}
\def\Srp{S^r_p\M}
\def\r{R_1}
\def\N{\mathcal N}
\def\tv{{tv}}
\def\th{{th}}
\def\h{\frak{h}}
\def\f{\frac{1}{\sqrt{r^2+\Vert K(\xi,\xi)\Vert^2}}}
\def\fR{\frak R}
\def\fRic{\frak{Ric}}
\def\Kxi{K(\xi,\xi)}
\def\Kxit{K(\xi,\xi)^{th}}
\def\Kxixi{\Vert\Kxi\Vert}
\def\Kxixit{\Vert\Kxi^{th}\Vert}

\title{The geometry of the tangent and sphere bundles over statistical manifolds}
%\thanks
\author{Barbara Opozda}
\let\runtitle\lala
\subjclass{ Primary: 53B05, 53B20, 53C05, 53C07, 53C20}

\keywords{statistical structure, Hessian structure, curvatures of
Riemannian and statistical structures, tangent bundle, sphere
bundle}

%\thanks{}

\address{Faculty of Mathematics and Computer Science UJ,
ul. prof. S. \L ojasiewicza  6, 30-348 Cracow, Poland}

\email{barbara.opozda@im.uj.edu.pl} \maketitle

\begin{abstract}
In the paper a Riemannian structure on the tangent bundle is defined
by using a statistical structure $(g,\nabla)$ on the base manifold.
Expressions for various curvatures  of the structure are derived.
Some  rigidity results of the structure are proved. The main goal of
the paper is to initiate the study of sphere bundles over
statistical manifolds. Basic formulas for the geometry are
established. Sphere bundles with small radii over compact manifolds
are studied.

\end{abstract}

\section{Introduction}
The geometry of the tangent bundle of a Riemannian manifold is an
old topic and has huge literature. The theory is of great importance
in mathematics and physics. It started in 1958 with a paper of
Sasaki (\cite{Sasaki}), who constructed a natural metric on the
tangent bundle $T\M$ by using a metric and its Levi-Civita
connection on the base manifold $\M$. The Levi-Civita connection
gives the horizontal bundle being a complementary distribution to
the canonically defined vertical bundle in $TT\M$. The horizontal
and vertical bundles are declared to be orthogonal relative to the
lifted metric. A horizontal distribution, however, can be defined by
means of any linear connection. This attempt was proposed by
Dombrowski in \cite{D} in 1962. He proved, in particular, that the
horizontal distribution is integrable if and only if the connection
on the base manifold is flat. Moreover, a naturally defined almost
complex structure on $TM$ is integrable if and only if the initial
connection is flat.

The metric  defined by Sasaki is  called the Sasaki metric and it is
a  fundamental concept in the theory. There exist, however,
important generalizations of this notion, e.g. the Cheeger-Gromoll
matric.

Since a statistical structure is a pair of a metric and a connection
(related in a way with the metric but very much independent of it),
it is natural to construct a metric on the tangent bundle by means
of the both objects.
%Recall that a statistical structure is a pair
%$(g,\nabla)$, where $g$ is a metric tensor field and $\nabla$ is a
%torsion-free connection such that the cubic form $\nabla g$ is
%symmetric.
The case where the connection is flat is very special in both
theories:  the theory initiated by Dombrowski and the one of
statistical structures. In the last theory the flat statistical
structures are Hessian structures, which are directly used in
statistical models. The structure on the tangent bundle determined
by a Hessian structure is Kaehlerian and its properties reflect
properties of the Hessian structure in a nice way. This case is
discussed in the book \cite{shima} and many papers.

In this paper we shall concentrate on the general case, that is, we
use an arbitrary statistical structure to construct (following the
idea of Sasaki) a metric on the tangent bundle. Some results in the
general case were proved in \cite{Mat} and \cite{Sat}. Note that the
general case is much more complicated than the Riemannan or the
Hessian one. For instance, one can  compare formulas for the
curvature tensor in respective cases, that is, formulas in
Proposition \ref{curvature_tR}, Corollary \ref{curvature_conj_symm}
for $K=0$ and in Corollary \ref{tR_Hessian}. In this paper we
establish basic formulas and facts for this theory. In particular,
we give explicit formulas for the curvature and Ricci tensors, for
the sectional  and scalar curvatures. Note that in the formulas
naturally  appear  some new curvatures  of statistical structures
introduced in \cite{BW4} and \cite{seccur}. From the formulas we
derive exemplary theorems generalizing the results known in the
Riemannian and Hessian cases. In particular, we prove that the
Riemannian structure on the tangent bundle is rigid in the sense
that if its sectional or scalar curvature is bounded, then the
structure must be flat see Theorems \ref{thm_scalar_curv_tan},
\ref{thm_sec_curv_tan}, \ref{thm_flatness_tan}.
%The Riemannian manifold $T\M$ has natural hypersurfaces,
%namely, sphere bundles $\Sr$ of various radii $r$, which are not so
%rigid.

The main aim of this paper is to initiate the study of sphere
bundles over statistical manifolds. To the knowledge of the author
of this paper, there are no papers dealing with this topic, even in
the Hessian case. Of course, there exist a lot of papers on the
sphere bundles in the Riemannian case. The complications, comparing
with the Riemannian case, are even bigger than those in the case of
tangent bundles. The complications start at the very beginning
because, unlike  the Riemannian  case, the horizontal distribution
is  not tangent to  a sphere bundle, in general, see Lemma
\ref{lemma_tangentspace_N}. In this paper we derive formulas for the
second fundamental and the mean curvature, (\ref{h}),
(\ref{mean_curvature}). We can also control  the scalar curvature,
although its expression is so long and complicated that, for the
purposes of this paper, it does not make much sense to write it down
explicitly. We consider sphere bundles with small radii. We prove
that if  $\M$ is compact, then for any given number $L\in \R$ there
exists  a sufficiently small radius $r$ such that the norm of the
mean curvature is greater than $L$ on the whole of $\Sr$, see
Theorem \ref{thm_normh}. If moreover $\dim\M>2$, then for a
sufficiently small radius  $r$ the scalar curvature of $\Sr$ is
greater than any prescribed number $L$ everywhere on $\Sr$, see
Theorem \ref{thm_scalar}.

\bigskip

\section{On the geometry of the tangent bundle determined by a statistical structure on the base manifold}
\subsection{Preliminaries on statistical manifolds} In this section
we collect basic information on statistical structures needed in
this paper. All proofs of the facts mentioned in this section can be
found in  \cite{BW4} and \cite{seccur}.

 A statistical structure on a manifold $\M$ is a pair $(g,\nabla)$,
where  $g$ is a metric tensor field (positive definite in this
paper) and $\nabla$ is a torsion-free connection such that the cubic
form $\nabla g$ is totally symmetric, that is, $(\nabla_X
g)(Y,Z)=(\nabla_Yg)(X,Z)$ for every $X,Y,Z$. Such a connection
$\nabla$ is called statistical . The Levi-Civita connection  for $g$
will be denoted by $\gnabla$. The difference tensor defined by
\begin{equation}\label{K_nabla}
K(X,Y)=K_XY=\nabla_XY-\gnabla_XY
\end{equation}
 is symmetric and symmetric relative to $g$. A statistical structure is called trivial
 if $\nabla=\gnabla$, that is, $K\equiv 0$. A statistical structure can be also defined as a pair
$(g,K)$, where $K$ is a $(1,2)$-tensor field symmetric and symmetric
relative to $g$. Namely, the corresponding statistical connection
$\nabla$ is defined by (\ref{K_nabla}). The cubic form $\nabla g$
and the difference tensor $K$ are related by the formula

\begin{equation}(\nabla g)(X,Y,Z)=-2g(K(X,Y),Z).
\end{equation}

 Having a metric tensor $g$ and a connection
$\nabla$ on a manifold $M$, one defines the conjugate connection
$\nablastar$ by the formula
\begin{equation}
g(\nabla _XY,Z)+g(Y,\nablastar _XZ)=Xg(Y,Z).
\end{equation}
A pair $(g,\nabla)$ is a statistical structure if and only if
$(g,\nablastar)$ is.

  The
curvature tensors for $\nabla$, $\nablastar$, $\gnabla$ will be
denoted by $R$, $ R^*$ and $R^g$, respectively. The corresponding
Ricci tensors will be denoted by $\Ric$, $\Ric^*$ and $\Ric^g$. A
statistical structure is Hessian if and only if $\nabla$ is flat,
i.e. $R=0$. For the fundamentals of the theory of Hessian structures
we refer to \cite{shima}. In general, the curvature tensor $R$ does
not satisfy the equality $g(R(X,Y)Z,W)=-g(R(X,Y)W,Z)$. If it does,
we say that the statistical structure is conjugate symmetric. The
notion was introduced in  \cite{L}. For any statistical structure we
have
\begin{equation}\label{RioR}
g(R(X,Y)Z,W)=-g(R^*(X,Y)W,Z).
\end{equation}
 The following
conditions are equivalent:
\newline
{\rm 1)} $R= R^*$,
\newline
{\rm 2)} $\gnabla K$ is symmetric, that is, $(\gnabla_X K)(Y,Z)$ is
symmetric for $X,Y$,
\newline
{\rm 3)} $g(R(X,Y)Z,W)$ is skew-symmetric relative to $Z,W$.

Therefore, each of the above three conditions characterizes a
conjugate symmetric statistical structure. The above equivalences
follow from the following well-known formula
\begin{equation}\label{from_Nomizu_Sasaki}
R(X,Y)= R^g(X,Y) +(\gnabla_XK)_Y-(\gnabla_YK)_X+[K_X,K_Y].
\end{equation}
Writing the same equality for $\nablastar$ and adding the both
equalities, we get
\begin{equation}\label{R+Rstar}
R+R^* =2R^g +2[K,K] ,\end{equation} where $[K,K]$ is a
$(1,3)$-tensor field defined as $[K,K](X,Y)Z=[K_X,K_Y]Z$.
 In particular, if $R=R^*$, then
\begin{equation}\label{z[K,K]dlahiperpowierzchni}
R=R^g +[K,K].
\end{equation}
Set
\begin{equation}
E=\tr_gK,\ \ \ \ \tau(X)=\tr K_X.
\end{equation}
The $1$-form $\tau$ is called the first Koszul form. The
$(0,2)$-tensor field $\nabla \tau$ is the second Koszul form.

For any connection $\nabla$ and any $(1,2)$-tensor field $K$ we have
\begin{equation}\label{about_traces}
\tr\nabla K(X,Y,\cdot)=\nabla\tau(X,Y),
\end{equation}
where $\tau(X)=\tr K_X$. Indeed, for any $(1,1)$-tensor field $S$
and any connection we have $X(\tr S)=\tr (\nabla _XS)$, see e.g.
\cite{BW4}. Applying this for $S=K_Y$ we get
\begin{equation}\label{1}
\nabla\tau(X,Y)=X(\tr K_Y)-\tr K_{\nabla_XY}=\tr (\nabla_X K_Y)-\tr
K_{\nabla_XY}.
\end{equation}
Let $e_1,...,e_n$ be a local frame around $p\in \M$. We have
\begin{equation}\label{2}
\begin{array}{rcl}
&&(\nabla _XK)(Y,e_i)=\nabla_X(K_Ye_i)-
K(\nabla_XY,e_i)-K(Y,\nabla_X e_i)\\
&&\ \ \ \ \ \ \ \ \ \ \ \ \ \ \ \ \ =(\nabla
_XK_Y)(e_i)-K(\nabla_XY,e_i).
\end{array}
\end{equation}
%Comparing (\ref{1}) and (\ref{2}) we get (\ref{about_traces}).
 We
shall use (\ref{about_traces}) for a statistical connection $\nabla$
as well as for the Levi-Civita connection $\gnabla$. Since both
connections are torsion-free, we have
$d\tau(X,Y)=\nabla\tau(X,Y)-\nabla\tau(Y,X)=\gnabla\tau(X,Y)-\gnabla\tau(Y,X)$.
By (\ref{from_Nomizu_Sasaki}) and (\ref{about_traces}) applied to
$\gnabla$ we have
\begin{equation}\label{dtau}
\tr
R(X,Y)=d\tau(X,Y)=\nabla\tau(X,Y)-\nabla\tau(Y,X)=\gnabla\tau(X,Y)-\gnabla\tau(Y,X).
\end{equation}
Since $\tr R(X,Y)=\Ric(X,Y)-\Ric(Y,X)$, the first Koszul form $\tau$
is closed if and only if the Ricci tensor $\Ric$ is symmetric. In
such a case the second Koszul form is symmetric. It happens, for
instance, if the structure is conjugate symmetric.

We also have the following equality
\begin{equation}\label{curvature_with_nablaK}
R(X,Y)=\gR(X,Y)+(\nabla_X K)_Y-(\nabla_Y K)_X-[K_X,K_Y].
\end{equation}
Using it and (\ref{about_traces}) we receive
\begin{equation}\label{for_Ric_vv}
\sum_{i=1}^n g(\nabla
K(e_i,Y,Z),e_i)=\Ric(Y,Z)-\gRic(Y,Z)+\RicK(Y,Z)+\nabla\tau(Y,Z).
\end{equation}

The $(1,3)$-tensor field $[K,K]$ has the same algebraic properties
as the Riemannian curvature tensor. In particular, $(g[K_X,K_Y]Z,W)$
is skew-symmetric for $Z,W$. Hence it can be used to define a
sectional curvature. Namely, if $\pi$ is a vector plane in $T_p\M$
and $e_1,e_2$ is an orthonormal basis of it, then
\begin{equation}\kK(\pi)=g([K,K](e_1,e_2)e_2,e_1).\end{equation}
This sectional curvature was studied in \cite{seccur}. We call this
curvature the sectional $K$-curvature. As in the Riemannian case, if
the sectional $K$-curvature vanish for all planes, then so does the
curvature tensor $[K,K]$ on $\M$. The curvature tensor $[K,K]$
determines the Ricci tensor, which will be denoted by $\Ric^K$. We
have
\begin{equation}
\Ric^K(X,Y)=\tau(K(X,Y))-g(K_X,K_Y).
\end{equation}

The curvature tensor
\begin{equation}
\frak R=\frac{R+R^*}{2}
\end{equation}
also has the same symmetries as the Riemannian curvature tensor and,
therefore,  it can be used to define a sectional curvature.
% In particular,
%\begin{equation}\label{frakR_property}
%g(\frak R(X,Y)W,Z)=g(\frak R(Z,W)Y,X).
%\end{equation}
  The Ricci tensor for
 $\frak R$ will be denoted by $\frak{Ric}$. Of course
 $2\frak{Ric}=\Ric+\Ric^*$. The Ricci tensors $\Ric^K$ and
 $\frak{Ric}$ are symmetric. The formula (\ref{R+Rstar})
 yields
 \begin{equation}\label{Ric+Ricstar}
2\frak{Ric}=\Ric+\Ric^*=2\Ric^g+2\Ric^K.
 \end{equation}
 Since $K$ is symmetric and symmetric relative to $g$,
 $g((\gnabla_XK)(Z,W),Y)$ is symmetric for $Z,W,Y$.
By (\ref{from_Nomizu_Sasaki}) we know that
\begin{equation*}
g(R(X,Y)\xi,\xi)=g((\gnabla_X K)(Y,\xi),\xi) -g((\gnabla_Y
K)(X,\xi),\xi).
\end{equation*}
Consequently,
\begin{equation}\label{for_sec_h}
\begin{array}{rcl}
&&\frac{1}{2}g(R(Y,X)\xi,\xi)+g((\gnabla_XK)(\xi,\xi),Y)\\
&&\ \ \ \ \ \ \ \ \ =\frac{1}{2}\{g((\gnabla_XK)(\xi,\xi),Y)+
g((\gnabla_YK)(\xi,\xi),X)\}.
\end{array}
\end{equation}

The scalar curvature $\rho$ of a statistical structure $(g,\nabla)$
is defined by  $\rho=\tr _g\Ric$. If we define the scalar curvature
for the dual structure $(g,\nabla ^*)$  we get the same $\rho$.
Hence
\begin{equation}
\rho=\tr_g\frak{Ric}.
\end{equation}
Of course we also have the scalar curvature $\rho ^g$ of $g$. We
have the following  relation between the scalar curvatures

\begin{equation}\rho^g=\rho+ \Vert K\Vert^2-\Vert
\tau\Vert^2.\end{equation}

\bigskip

\subsection{The Sasaki metric on the tangent bundle induced by a statistical structure}
 Let $\M$ be an $n$-dimensional manifold. Consider the tangent
bundle $\pi:T\M\to \M$. $T\M$ is a $2n$-dimensional manifold. Denote
by $TT\M$ the tangent bundle of $T\M$ and by $\pi_*$ the
differential  $\pi_*:TT\M\to T\M$. For each point $\xi\in \T\M$ the
space $\mathcal V_{\xi}=\ker (\pi_*)_\xi$ is called the vertical
space at $\xi$. It  is a $n$-dimensional vector subspace of $T_\xi
T\M$. The space $\mathcal V_\xi$ can  be naturally
 identified with $T_{\pi(\xi)}\M$. Namely,  let $p=\pi (\xi)$ and
 $X\in T_p\M$. We define the  vertical  lift $X^v_\xi$ of $X$ to $\xi$ by saying that
 the curve $t\to \xi+tX$ lying in the affine space $T_p\M$ is an integral
 curve of the vector $X^v_\xi$. We shall write it  as follows:
 $X^v_\xi=[t\to \xi+tX]$.
It is clear that $\mathcal V_\xi=\{ X^v_\xi: \ \ X\in T_p\M \}$ and
the assignment $T_p\M\ni X\to X^v_\xi\in  \mathcal V_\xi$ is a
linear isomorphism. The vertical bundle $\mathcal V=\bigcup_{\xi \in
T\M}\mathcal V_\xi$ is a smooth $n$-dimensional vector subbundle of
$TT\M$, so it is a distribution in $TT\M$. It is clearly integrable.

Let $\nabla$ be a connection on $\M$. The connection  determines the
horizontal distribution  $\mathcal H=\bigcup _{\xi\in T\M}\mathcal
H_\xi $ in $TT\M$ complementary to $\mathcal V$, that is, $\mathcal
H_\xi\oplus\mathcal V_\xi=T_\xi T\M$ for  each $\xi\in T\M$. We
briefly explain how to obtain the horizontal  distribution. For
details we refer, e.g. to \cite{D}. One can produce the horizontal
lift $X^h_\xi$ of a vector $X\in T_p\M$ as follows. Take a geodesic
(relative to $\nabla$) $\gamma(t)$ determined by $X$. By using the
parallel transport of the vector $\xi$ along $\gamma$ relative to
$\nabla$ we obtain a curve $t\to \xi_t$ in $T\M$. Of course
$\xi_0=\xi$. We set $X^h_\xi=[t\to\xi_t]$. We see that $ (\pi_*)_\xi
(X^h_\xi)=[t\to\pi(\xi_t)]=[t\to \gamma(t)]=X$. The horizontal space
at $\xi$ is defined as $\mathcal H_\xi=\{X^h_\xi ;\ \ X\in T_p\M\}$.
The assignment $X\in T_p\M\to X^h_\xi\in \mathcal H_\xi$ is a linear
isomorphism. If $X$ is a vector field on $\M$, then $X^v$ and $X^h$
will stand for the vertical  and horizontal lifts to $T\M$. The
vertical and horizontal lifts are  smooth vector fields on $T\M$. Of
course, $X$ can be a locally defined vector field, say on $\mathcal
U\subset \M$, and then the lifts are defined on $T\M_{|\mathcal U}$.

For a smooth function $f$ on $\M$  we define the function  $F=f\circ
\pi$ on $T\M$. It is clear that
\begin{equation}\label{lifts_with_function}
F\circ X^a=(fX)^a,
\end{equation}
where $a=v$ or $h$, and
\begin{equation}\label{F=f_circ_pi}
 X^vF=0,\ \ \ \  \ (X^hF)_\xi=(Xf)_{\pi(\xi)}.
\end{equation}
We have the following  description of  the Lie bracket of vertical
and horizontal lifts, see \cite{D},
\begin{lemma}\label{Dabrowski1}
If $X,Y\in\mathcal X(\M)$, then  we have
\begin{equation}
[X^v,Y^v]=0,\ \ \ \ [X^h,Y^v]=(\nabla_XY)^v,\ \ \ \ \
[X^h,Y^h]_\xi=[X,Y]^h_\xi-(R(X,Y)\xi)^v_\xi ,
\end{equation}
where $R$ is the curvature tensor of $\nabla$.
\end{lemma}
The connection $\nabla$ determines also an almost complex structure
$J$ on the manifold $T\M$ as follows
\begin{equation}
JX^v=X^h,\ \ \ \ \ \  JX^h=-X^v  \ \ \ \ for \ \ \ \ X\in \mathcal
X(\M).
\end{equation}
Using Lemma \ref{Dabrowski1} one gets, see \cite{D},
\begin{proposition}\label{Dombrowski2}
The almost complex structure $J$ is integrable if and only if the
connection $\nabla$ is without torsion and its curvature tensor $R$
vanishes.
\end{proposition}
\medskip

Assume now that  we have a statistical structure $\stat$ on $\M$.
The metric tensor field $g$ can be lifted to $T\M$ in various ways.
We shall follow the easiest way proposed by Sasaki, \cite{Sasaki}.
We shall use both $g$ and $\nabla$ for constructing the lift.
Namely, we declare that the vertical and horizontal distributions
 are orthogonal everywhere  on $T\M$. Moreover, we declare that  the horizontal and
vertical lifts are both isometries. Denote the obtained metric by
$\tg$. Hence we have
\begin{equation}
\tg (X^v,Y^v)=g(X,Y)=\tg (X^h,Y^h),\ \ \ \ \  \tg (X^v,Y^h)=0.
\end{equation}
We call this metric the Sasakian metric. It is clear that $\tilde g$
is almost Hermitian relative to $J$. Hence we have the  almost
symplectic form $\omega$ on $T\M$ determined by the almost Hermitian
structure $(g,J)$. It is given by
\begin{equation*}
\omega(X^v,Y^v)=0, \ \ \ \ \omega(X^h,Y^h)=0,\ \ \ \ \  \
\omega(X^v,Y^h)=-\omega(Y^h,X^v)=g(X,Y).
\end{equation*}
Using Lemma \ref{Dabrowski1} and the first Bianchi identity for the
dual connection one easily sees that for any statistical structure
the form $\omega$ is closed, that is, the almost Hermitian structure
$(g,J)$ is almost Kaehlerian, see \cite{Sat}. Indeed, it is
sufficient to  check that  $d\omega(X^v,Y^v,Z^v)=0,\
d\omega(X^v,Y^v,Z^h)=0,\ d\omega(X^h,Y^h,Z^v)=0,\
d\omega(X^h,Y^h,Z^h)=0. $
 The first three
equalities immediately follow from Lemma \ref{Dabrowski1}. As
regards the fourth one we have
\begin{eqnarray*}
&&d\omega(X^h,Y^h, Z^h)=X^h\omega (Y^h,
Z^h]-Y^h\omega(X^h,X^h)+Z^h\omega(X^h,Y^h)\\
&&\ \ \ \ \ \ \ \ \ \ \ -\omega([X^h,Y^h],Z^h)+\omega([X^h,Z^h],Y^h)-\omega([Y^h,Z^h],X^h)\\
&&\ \ \ \ \ \ \ \ \ =-g(R(X,Y)\xi,Z)-g(R(Z,X)\xi,Y)-g(R(Y,Z)\xi,X)\\
&&\ \ \ \ \ \ \ \ \
=g(R^*(X,Y)Z,\xi)+g(R^*(Z,X)Y,\xi)+g(R^*(Y,Z)X,\xi)=0.
\end{eqnarray*}
By Proposition \ref{Dombrowski2} it is now clear  that the almost
Hermitian structure on the tangent bundle is Kaehlerian if and only
if the given statistical structure is Hessian.

\medskip

\medskip

Denote by $\tnabla$ the Levi-Civita connection of $\tilde g$. Using
the Koszul formula for the covariant derivative of the Levi-Civita
connection, Lemma \ref{Dabrowski1} and (\ref{F=f_circ_pi}) we get

\begin{lemma}\label{tnabla}
The Levi-Civita connection of $\tg$ is given by the following
formulas
\begin{equation}
\begin{array}{rcl}
&&(\tnabla_{X^h}Y^h)_\xi=-\frac{1}{2}(R(X,Y)\xi)^v_\xi+(\gnabla_XY)^h_\xi\\
&&\ \ \ \ \ \ \ \ \ \ \ \ \ \ \ \ \  \ =
-\frac{1}{2}(R(X,Y)\xi)^v_\xi+(\nabla_XY)^h_\xi- K(X,Y)^h_\xi,
\end{array}
\end{equation}
\begin{equation}
(\tnabla_{X^h}Y^v)_\xi=(\gnabla_XY)^v_\xi-\r(X,\xi,Y)^h_\xi=(\nabla_XY)^v_\xi-K(X,Y)^v_\xi-\r(X,\xi,Y)^h_\xi,
\end{equation}

\begin{equation}
(\tnabla_{Y^v}X^h)_\xi=-K(X,Y)^v_\xi-\r(X,\xi,Y)^h_\xi,
\end{equation}
\begin{equation}
(\tnabla_{X^v} Y^v)_\xi=K(X,Y)^h_\xi,
\end{equation}
where $X,Y\in\mathcal X(\M)$ and $\r(X,\xi,Y)$ is defined by the
formula $g(\r(X,\xi,Y),Z)=\frac{1}{2} g(R(Z,X)\xi,Y)$ for every $Z$.
\end{lemma}

For technical reasons we have introduced the  tensor field $\r$ of
type $(1,3)$ defined by
\begin{equation}
g(W,\r(X,Y,Z))=\frac{1}{2}g(R(W,X)Y,Z).
\end{equation}
If a  statistical structure is conjugate symmetric, then
$g(R(W,X)Y,Z)=g(R(Y,Z) W,X)$ and thus
\begin{equation}
\r(X,Y,Z)=\frac{1}{2}R(Z,Y)X.
\end{equation}
We shall need the following equalities holding for any statistical
structure

\begin{equation}\label{for_curvature_tg}
\r(Z,Y,X)-\r(Z,X,Y)=\frak R(X,Y)Z=R^g(X,Y)Z+[K_X,K_Y]Z.
\end{equation}
\bigskip
To check (\ref{for_curvature_tg}), it is sufficient to observe
(using the symmetries of $\frak R$ and (\ref{R+Rstar})) that
$$g(\r(Z,Y,X)-\r(Z,X,Y), W)=g(\frak R(W,Z)Y,X)=g(\frak R(X,Y)
Z,W).$$

\bigskip

\subsection{Curvatures of the Sasakian metric determined by a statistical
structure} Let now $\F$ be a $(1,1)$-tensor field on $\M$ (of
course, maybe defined on some open subset of $\M$). Let $X,\xi\in
T_p\M$ and $\tilde Y$ be a vector field on $T\M$ defined at least
along an integral curve of $X^a_\xi$, where $a$ stands for $v$ or
$h$.
%We shall now look  at $(\tnabla_{X^a}\F(Y)^b)_\xi$, where $\xi\in T_p\M$ and $a,b$ stand for $v$ or $h$.
If $X^a_\xi=[t\to \gamma (t)]$, then $(\tnabla_{X^a}\tilde
Y)_\xi=\tnabla_{X^a_\xi}\tilde Y_{\gamma ( t)}$.
% We want to
%compute
%$(\tnabla_{X^a}\F(Y)^b)_\xi=\tnabla_{X^a_\xi}\F(Y)^b_{\xi_t}$.
Assume that $\tilde Y_{\gamma(t)}=((\F\circ \gamma)(t))^b$, where
$b$ denotes  $v$ or $h$. Let $a=v$. Then $X^v_\xi=[t\to
\gamma(t)=\xi+tX]$. We have $\F(\gamma(t))=\F(\xi)+t\F(X)\in T_p\M$,
and consequently, by Lemma \ref{tnabla} we get
\begin{eqnarray*}
\tnabla_{X^v_\xi}((\F\circ\gamma)(t))^v=\tnabla_{X^v_\xi}(\F(\xi)^v+t\F(X)^v)=K(X,\F(\xi))^h_\xi+\F
(X)^v_\xi,\end{eqnarray*}
\begin{eqnarray*}
&&\tnabla_{X^v_\xi}((\F\circ\gamma)(t))^h
=\tnabla_{X^v_\xi}(\F(\xi)^h+t\F(X)^h)\\
&&\ \ \ \ \ \ \ \ \ \ \ \ \ \ \
=-K(\F(\xi),X)^v_\xi-\r(\F(\xi),\xi,X)^h_{\xi}+\F(X)^h_\xi.
\end{eqnarray*}
 We now consider the case, where $a=h$. Extend  the vector $\xi$
along geodesics (relative to $\nabla$) by the parallel transport. We
denote the obtained vector field (on a domain of a normal coordinate
system) by  $\Xi$. We have $\Xi(p)=\xi$ and $(\nabla\Xi)_p=0$. Let $
t\to \rho(t)$ be the geodesic for $X$. Then
$X^h_\xi=[t\to\gamma(t)=\Xi(\rho(t))]$, and consequently,
$(\tnabla_{X^h_\xi}(\F\circ\gamma)
^b=\tnabla_{X^h_\xi}(\F(\Xi(\rho(t)))^b$. We have

\begin{eqnarray*}
&&\tnabla_{X^h_\xi}(\F(\Xi(\rho(t)))^h=
-\frac{1}{2}(R(X,\F(\xi))\xi)^v_\xi+(\nabla_X\F(\Xi))^h_\xi-K(X,\F(\Xi))^h_\xi\\
&&\ \ \ \ \ \ \ \ \ \ \ \ \ \ \ \ \ \ \ \
=-\frac{1}{2}(R(X,\F(\xi))\xi)^v_\xi+((\nabla_X\F)(\xi))^h_\xi-K(X,\F(\xi))^h_\xi,
\end{eqnarray*}

\begin{eqnarray*}
&&\tnabla_{X^h_\xi}(\F(\Xi(\rho(t)))^v=
-\r(X,\xi,\F(\xi))^h_\xi+(\nabla_X\F(\Xi))^v_\xi-K(X,\F(\xi))^v_\xi\\
&&\ \ \ \ \ \ \ \ \ \ \ \ \ \ \ \ \ \ \ \
=-\r(X,\xi,\F(\xi))^h_\xi+((\nabla_X\F)(\xi))^v_\xi-K(X,\F(\xi))^v_\xi.
\end{eqnarray*}
\begin{lemma}\label{tensor_of_type11}
Let $\F$ be a $(1,1)$-tensor field on $\M$ and $X,\xi\in T_p\M$,
$p\in M$. If $\gamma (t)=\xi+tX$, then
\begin{equation}
\tnabla_{X^v_\xi}((\F\circ\gamma)(t))^v=K(X,\F(\xi))^h_\xi+\F
(X)^v_\xi,
\end{equation}
\begin{equation}
\tnabla_{X^v_\xi}((\F\circ \gamma
)(t))^h=-K(\F(\xi),X)^v_\xi-\r(\F(\xi),\xi,X)^h_{\xi}+\F(X)^h_\xi.
\end{equation}
If $\gamma(t)=\Xi(\rho(t))$, where $\rho(t)$ is a geodesic
determined by $X$, then
\begin{equation}
\tnabla_{X^h_\xi}((\F\circ\gamma)(t))^h=-\frac{1}{2}(R(X,\F(\xi))\xi)^v_\xi+((\nabla_X\F)(\xi))^h_\xi-K(X,\F(\xi))^h_\xi,
\end{equation}
\begin{equation}
\tnabla_{X^h_\xi}((\F\circ \gamma )(t))^v
=-\r(X,\xi,\F(\xi))^h_\xi+((\nabla_X\F)(\xi))^v_\xi-K(X,\F(\xi))^v_\xi.
\end{equation}
\end{lemma}

For later use observe that in the same manner as above one can prove

\begin{lemma}\label{tensor_of_type1k}
Let $\F$ be a $(1,k)$-tensor field on $\M$ and $X,\xi\in T_p\M$,
$p\in M$. If $\gamma (t)=\xi+tX$, then
\begin{equation}
\begin{array}{rcl}
&&\tnabla_{X^v_\xi}(\F(\gamma(t),...,\gamma(t)))^v=K(X,\F(\xi,...,
\xi))^h_\xi\\
&&\ \ \ \ \ \ \ \ \ \ \ \ \ \ \ +[\F (X,\xi,
...,\xi)+\F(\xi,X,...,\xi)+\F(\xi,...,\xi,X)]^v_\xi,
\end{array}
\end{equation}
\begin{equation}
\begin{array}{rcl}
&&\tnabla_{X^v_\xi} (\F(\gamma(t),...,\gamma(t)))^h=-K(X,\F(\xi,...,
\xi))^v_\xi-\r(\F(\xi,...,\xi),\xi,X)^h_\xi\\
&&\ \ \ \ \ \ \  \ \ \ \ \ \ \ \ +[\F (X,\xi,
...,\xi)+\F(\xi,X,...,\xi)+\F(\xi,...,\xi,X)]^h_\xi.
\end{array}
\end{equation}
If $\gamma(t)=\Xi(\rho(t))$, where $\rho(t)$ is a geodesic
determined by $X$, then
\begin{equation}
\begin{array}{rcl}
&&\tnabla_{X^h_\xi}(\F (\gamma(t),...,\gamma(t)))^h=-\frac{1}{2}(R(X,\F(\xi,...,\xi))\xi)^v_\xi\\
&&\ \ \ \ \ \ \ \ \ \ \ \
+((\nabla_X\F)(\xi,...,\xi))^h_\xi-K(X,\F(\xi,...,\xi))^h_\xi,
\end{array}
\end{equation}
\begin{equation}
\begin{array}{rcl}
&&\tnabla_{X^h_\xi}(\F(\gamma(t),...,\gamma(t)))^v
=-\r(X,\xi,\F(\xi,...,\xi))^h_\xi\\
&&\ \ \ \ \ \ \ \ \
+((\nabla_X\F)(\xi,...,\xi))^v_\xi-K(X,\F(\xi,...,\xi))^v_\xi.
\end{array}
\end{equation}
\end{lemma}
In particular, if $\F=K$, using the above lemma and the equality
$\nabla=\gnabla+K$, we get
\begin{equation}\label{tensor_K_vv}
\tnabla_{X^v_\xi}(K(\gamma(t),\gamma(t)))^v=K(X,K(\xi, \xi))^h_\xi+2
K(X,\xi)^v_\xi,
\end{equation}

\begin{equation}\label{tensor_K_vh}
\tnabla_{X^v_\xi}
(K(\gamma(t),\gamma(t)))^h=-K(X,K(\xi,\xi))^v_\xi-\r(K(\xi,\xi),\xi,X)^h_\xi
+2K(X,\xi)^h_\xi,
\end{equation}

\begin{equation}\label{tensor_K_hh}
\tnabla_{X^h_\xi}(K (\gamma(t),\gamma(t)))^h=-\frac{1}{2}(R(X,
K(\xi,\xi))\xi)^v_\xi+ ((\gnabla_X
K)(\xi,\xi))^h_\xi-2(K_\xi^2(X))^h_\xi,
\end{equation}

\begin{equation}\label{tensor_K_hv}
\tnabla_{X^h_\xi}(K(\gamma(t),\gamma(t)))^v
=-\r(X,\xi,K(\xi,\xi))^h_\xi +((\gnabla_X
K)(\xi,\xi))^v_\xi-2(K_\xi^2(X))^v_\xi.
\end{equation}

\bigskip

Using  Lemmas \ref{Dabrowski1}, \ref{tnabla}, \ref{tensor_of_type11}
as well as  the properties of statistical structures (collected in
Section 2.1) and affine connections (in particular, the second
Bianchi identity for $\nabla$), we can now compute the curvature
tensor $\tR$ for $\tg$. In the following computations we adopt the
following conventions. We fix $\xi, X,Y,Z \in T_p\M$. Extend them
along $\nabla$-geodesics by the parallel transport to vector fields
(denoted by the same letters) defined in a normal neighborhood of
$p$. In particular, $[X,Y]_p=0$, $\nabla_{X_p}Y=0$ and
$\gnabla_{X_p}Y=-K(X,Y)_p$. The vector fields are now fixed. If we
have a $(1,1)$-tensor field $\mathcal F$, then the expression
$\tnabla_{{X^a}_\xi}(\mathcal F(\gamma))^b$ stands for
$\tnabla_{{X^a}_\xi}(\mathcal F\circ\gamma)^b$, where $\gamma$ is an
integral curve of $X^a_\xi$ defined as above and where $a,b$ are $v$
or $h$. In the following computations we use the following
$(1,1)$-tensor fields: $\mathcal F(U)= \r(Z,U, Y)$, $\mathcal
F(U)=\r(Y,U,Z)$, $\mathcal F(U)=\r(Z,U,X)$ or $\mathcal
F(U)=R(Y,Z)U$. We shall omit $\xi$ in the subscript position. We
have

\begin{equation}
\begin{array}{rcl}
&&\tnabla_\vX\tnabla_\vY\vZ=\tnabla_\vX
(K(Y,Z))^h\\&&=-(\r(K(Y,Z),\u,X))^h-(K_X(K(Y,Z)))^v,
\end{array}
\end{equation}

\begin{equation}
\begin{array}{rcl}
&&\tnabla_\vX\tnabla_\vY\hZ=\tnabla_\vX\left(-\r(Z,\gamma,Y)^h-K(Z,Y)^v\right)\\
&&
=K(\r(Z,\u,Y),X)^v+\r(\r(Z,\u,Y),\u,X)^h-\r(Z,X,Y)^h\\
&&\ \ \ \ -(K_X(K(Y,Z)))^h,
\end{array}
\end{equation}

\begin{equation}
\begin{array}{rcl}
&&\tnabla_{\vX}\tnabla_\hY\hZ=\tnabla_{\vX}\left(-\frac{1}{2}(R(Y,Z)\gamma)^v+(\nabla_YZ)^h-K(Y,Z)^h\right)\\
&&=-\frac{1}{2} \left(
(K(X,R(Y,Z)\u))^h+(R(Y,Z)X)^v\right)\\
&&+r(K(Y,Z),\u,X)^h+K_X(K(Y,Z))^v,
\end{array}
\end{equation}

\begin{equation}\label{for_exemplary_computation}
\begin{array}{rcl}
&&\tnabla_\hY\tnabla_\vX\hZ=\tnabla_\hY\left(-\r(Z,\gamma,X)^h-K(X,Z)^v\right)\\
&&
=\frac{1}{2}(R(Y,\r(Z,\u,X))u)^v-((\nabla_Y\r)(Z,\u,X))^h\\
&& \ \ \ \ \ \ \ \ +(K(Y,\r(Z,\u,X)))^h -\nabla K
(Y,X,Z)^v\\
&&\ \ \ \ \ \ \ \ +(K_YK(X,Z))^v+(\r(Y,\u,K(X,Z))^h,
\end{array}
\end{equation}

\begin{equation}
\begin{array}{rcl}
&&\tnabla_\vX\tnabla_\hY\vZ=\tnabla_\vX(\nabla_YZ-K(Y,Z))^v-\tnabla
_\vX \r(Y,\gamma,Z)^h\\
&&=-(K_XK(Y,Z)))^h\\
&&-\left( -K(\r(Y,\u,Z),X)^v-\r(\r(Y,\u,Z),\u,X)^h+\r(Y,X,Z)^h
\right),
\end{array}
\end{equation}

\medskip

\begin{equation}
\begin{array}{rcl}
&&\tnabla_\hY\tnabla_\vX\vZ=\tnabla_\hY (K(X,Z))^h\\
&& =-\frac{1}{2} (R(Y,K(X,Z))\u)^v+((\nabla
K)(Y,X,Z))^h-(K_Y(K(X,Z)))^h,
\end{array}
\end{equation}

\medskip

\begin{equation}
\begin{array}{rcl}
&&\tnabla_{\hX}\tnabla_{\hY}\hZ= \tnabla_\hX\left(
-\frac{1}{2}(R(Y,Z)\gamma)\u^v+(\gnabla_YZ)^h\right)\\
&&
=-\frac{1}{2}\left(-\r(X,\u,R(Y,Z)\u)^h+((\nabla_XR)Y,Z)\u)^v-K(X,R(Y,Z)\u)^v\right)\\
&&-\frac{1}{2}(R(X,\gnabla_YZ)\u)^v+\left(\gnabla_X\gnabla_YZ
\right)^h,
\end{array}
\end{equation}

\begin{equation}
\begin{array}{rcl}
&&\tnabla_{[\hX,\hY]}\hZ=\tnabla_{[X,Y]^h}Z^h-\tnabla_{(R(X,Y)\u)^v}\hZ\\
&&=\r(Z,\u,R(X,Y)\u))^h+K(Z,R(X,Y)\u)^v,
\end{array}
\end{equation}

\medskip

\begin{equation}
\begin{array}{rcl}
&&\tnabla_\hX\tnabla_\hY\vZ=\tnabla_\hX\left((\gnabla_YZ)^v-\r(Y,\gamma,Z)^h\right)\\
&&  =(\gnabla_X\gnabla_YZ)^v+\r(X,\u,K(Y,Z))^h\\
&& - [- \frac{1}{2} (R(X,\r(Y,\u,Z))\u)^v +\nabla_X\r)(Y,\u,Z)^h\\
&&\ \ \ \ \ \ \ \ \ \ \ \ \ \ \ \ \ \ \ \ \ \ \ \ \ \
-K(X,\r(Y,\u,Z))^h],
\end{array}
\end{equation}
\medskip

\begin{equation}
\tnabla_{[\hX,\hY]}\vZ=-\tnabla_{(R(X,Y)\u)^v}\vZ=-K(R(X,Y)\u,Z)^h.
\end{equation}

From the above formulas one easily gets the following expressions
for the curvature tensor $\tR$. Note that  for proving (\ref{vvh})
one also needs (\ref{for_curvature_tg}). For proving (\ref{hhh}) the
second Bianchi identity for $\nabla$ is needed.

\begin{proposition}\label{curvature_tR}
The curvature tensor $\tR_\u$ of $\tg$ is given by the following
formulas
\medskip
%1
\begin{equation}
\begin{array}{rcl}
&&\tR(\vX,\vY)\vZ=([K_Y,K_X]Z)^v\\
&&\ \ \ \ \ \ \ \ \ \ \ \ \ \ \
+\{\r(K(X,Z),\u,Y)-\r(K(Y,Z),\u,X)\}^h,
\end{array}
\end{equation}
\medskip

%2
%\begin{equation}
%\begin{array}{rcl}
%&&\tR(\vX,\vY)\hZ=\{K(r(Z,\u,Y),X)-K(r(Z,\u,X),Y)\}^v\\
%&&\ \ \ \ \ \ \ \ \ \ \ \ +\{[K_Y,K_X]Z+ r(r(Z,\u,Y),\u,X)-r(r(Z,\u,X),\u,Y)\\
%&&\ \ \ \ \ \ \ \ \ \ \ \ \ \ \ \ \ \ \ \ \ \ \ \ \
%+r(Z,Y,X)-r(Z,X,Y)\}^h,
%\end{array}
%\end{equation}
%\medskip

%2
\begin{equation}\label{vvh}
\begin{array}{rcl}
&&\tR(\vX,\vY)\hZ=\{K(\r(Z,\u,Y),X)-K(\r(Z,\u,X),Y)\}^v\\
&&\ \ \ \  +\{R^g(X,Y)Z+
\r(\r(Z,\u,Y),\u,X)-\r(\r(Z,\u,X),\u,Y)\}^h,
\end{array}
\end{equation}
\medskip

%3
\begin{equation}
\begin{array}{rcl}
&&\tR(\vX,\hY)\hZ=\{-\frac{1}{2}R(Y,Z)X+[K_X,K_Y]Z+\nabla
K(Y,X,Z)\\
&&\ \ \ \ \ \ \ \ \ \ \ \ \ \ \ \ \ \ \ \ \ \ \ \ \ \ -\frac{1}{2}R(Y,\r(Z,\u,X))\u\}^v\\
&&\ \ \ \ \ \ +\{-\frac{1}{2}K(X,R(Y,Z)\u)+\r(K(Y,Z),\u,X)+(\nabla_Y
\r)(Z,\u,X)\\
&&\ \ \ \ \ \ \ \ \ \ \ \ \ \ \ \ \ \
-K(Y,\r(Z,\u,X))-\r(Y,\u,K(X,Z))\}^h,
\end{array}
\end{equation}
\medskip

%4
\begin{equation}
\begin{array}{rcl}
&&\tR(\vX,\hY)\vZ=\{K(\r(Y,\xi,Z),X)+\frac{1}{2}R(Y,K(X,Z))\u\}^v\\
&&\ \ \ \ \ \ \ \ \ \ \ \ \  +\{ [K_Y,K_X]Z
 +\r(\r(Y,\u,Z),\u,X)\\
 &&\ \ \ \ \ \ \ \ \ \ \ \ \ \ \ \ \ \ \ \ \ \ \  \
 -\r(Y,X,Z)-(\nabla K)(Y,X,Z)\}^h,
\end{array}
\end{equation}
\medskip

%5
\begin{equation}
\begin{array}{rcl}
&&\tR(\hX,\hY)\hZ=(\gR(X,Y)Z)^h\\
&&
+\frac{1}{2}\{\r(X,\u,R(Y,Z)\u)-\r(Y,\u,R(X,Z)\u)-2\r(Z,\u,R(X,Y)\u)\}^h\\
&&\ \ \ \ \
+\frac{1}{2}\{(\nabla_ZR)(X,Y)\u-R(Y,K(X,Z))\u+R(X,K(Y,Z))\u\\
&&\ \ \ \ \ \ \ \ \ \
+K(X,R(Y,Z)\u)-K(Y,R(X,Z)\u)-2K(Z,R(X,Y)\u)\}^v,
\end{array}
\end{equation}
\medskip

%6
\begin{equation}
\begin{array}{rcl}
&&\tR(\hX,\hY)\vZ=(\gR(X,Y)Z)^v\\
&&\ \ \ \ \ \ \ +\frac{1}{2}\{R(X,\r(Y,\u,Z))\u-R(Y,\r(X,\u,Z))\u\}^v\\
&&\ \ \ \
+\{(\nabla_Y\r)(X,\u,Z)-(\nabla_X\r)(Y,\u,Z)\\
&&\ \ \ \ \ \ \ \ \ +K(X,\r(Y,\u,Z))-K(Y,\r(X,\u,Z)
+K(R(X,Y)\u,Z)\\
&&\ \ \ \ \ \ \ \ \ \ +\r(X,\u,K(Y,Z))-\r(Y,\u,K(X,Z))\}^h.
\end{array}
\end{equation}

\end{proposition}

\begin{corollary}\label{curvature_conj_symm}
 If $\stat$ is conjugate symmetric, then the
curvature tensor $\tR_\u$ is given by the following formulas:
%11

\begin{equation}
\begin{array}{rcl}
&&\tR(\vX,\vY)\vZ=([K_Y,K_X]Z)^v\\
&&\ \ \ \ \ \ \ \ \ \ \ \ \ \ \ +\{\frac{1}{2}
R(Y,\u)K(X,Z)-\frac{1}{2}R((X,\u)K(Y,Z)\}^h,
\end{array}
\end{equation}

\medskip
%12
\begin{equation}
\begin{array}{rcl}
&&\tR(\vX,\vY)\hZ=\{\frac{1}{2}K(R(Y,\u)Z,X)-\frac{1}{2}K(R(X,\u)Z,Y)\}^v\\
&&\ \ \ \ \ \ \ \  +\{\gR(X,Y)Z+
\frac{1}{4}R(X,\u)(R(Y,\u)Z)-\frac{1}{4}R(Y,\u)(R(X,\u)Z)\}^h,
\end{array}
\end{equation}

\medskip
%13
\begin{equation}
\begin{array}{rcl}
&&\tR(\vX,\hY)\hZ=\{-\frac{1}{2}R(Y,Z)X+[K_X,K_Y]Z+\nabla
K(Y,X,Z)\\
&&\ \ \ \ \ \ \ \ \ \ \ \ \ \ \ \ \ \ \ \ \ \ \ \ \ \ +\frac{1}{4}R(R(X,\u)Z,Y)\u\}^v\\
&&\ \ \ \ \ \
+\{-\frac{1}{2}K(X,R(Y,Z)\u)+\frac{1}{2}R(X,\u)K(Y,Z)+\frac{1}{2}(\nabla_YR)
)(X,\u)Z\\
&&\ \ \ \ \ \ \ \ \ \ \ \ \ \ \ \ \ \
-\frac{1}{2}K(Y,R(X,\u)Z))-\frac{1}{2}R(K(X,Z),\u)Y\}^h,
\end{array}
\end{equation}

\medskip
%14
\begin{equation}
\begin{array}{rcl}
&&\tR(\vX,\hY)\vZ=\{\frac{1}{2}K(R(Z,\u)Y,X)+\frac{1}{2}R(Y,K(X,Z))\u\}^v\\
&&\ \ \ \ \ \ \ \ \ \ \ \ \  +\{ [K_Y,K_X]Z
 +\frac{1}{4}R(X,\u)(R(Z,\u)Y)\\
 &&\ \ \ \ \ \ \ \ \ \ \ \ \ \ \ \ \ \ \ \ \ \ \  -\frac{1}{2}R(Z,X)Y-(\nabla K)(Y,X,Z)\}^h,
\end{array}
\end{equation}
\medskip
%15
\begin{equation}\label{hhh}
\begin{array}{rcl}
&&\tR(\hX,\hY)\hZ=\{(\gR(X,Y)Z)^h\\
&&\ \ \ \
+\frac{1}{4}R(R(Y,Z)\u,\u)X-\frac{1}{4}R(R(X,Z)\u,\u)Y-\frac{1}{2}R(R(X,Y)\u,\u)Z\}^h\\
&&\ \ \ \ \
\frac{1}{2}\{(\nabla_ZR)(X,Y)\u-R(Y,K(X,Z))\u+R(X,K(Y,Z))\u\\
&&\ \ \ \ \ \ \ \ \ \
+K(X,R(Y,Z)\u)-K(Y,R(X,Z)\u)-2K(Z,R(X,Y)\u)\}^v,
\end{array}
\end{equation}
\medskip

%16
\begin{equation}
\begin{array}{rcl}
&&\tR(\hX,\hY)\vZ=(\gR(X,Y)Z)^v\\
&&\ \ \ \ \ \ \ +\frac{1}{4}\{ R(X,R(Z,\u)Y)\u-R(Y,R(Z,\u)X)\u\}^v\\
&&\ \ \ \
+\frac{1}{2}\{(\nabla_YR)(Z,\u)X-(\nabla_XR)(Z,\u)Y\\
&&\ \ \ \ \ \ \ \ \ +K(X,R(Z,\u)Y)-K(Y,R(Z,\u)X)
+2K(R(X,Y)\u,Z)\\
&&\ \ \ \ \ \ \ \ \ \ \ \ \ \ \ \ \ \ \ \  \ \ \ \ \ \ \ \ \ \ \
+R((K(Y,Z)\u,X)\}^h.
\end{array}
\end{equation}
\end{corollary}

\begin{corollary}\label{tR_at_0}
The curvature tensor $\tR$ at $\u=0_p$ is given by
%21
\begin{equation}
\tR(\vX,\vY)\vZ=([K_Y,K_X]Z)^v,
\end{equation}
%22
\begin{equation}
\tR(\vX,\vY)\hZ= \{[K_Y,K_X]Z +\r(Z,Y,X)-\r(Z,X,Y)\}^h,
\end{equation}
%23
\begin{equation}
\tR(\vX,\hY)\hZ=\{-\frac{1}{2}R(Y,Z)X+[K_X,K_Y]Z+\nabla K(Y,X,Z)\}^v
\end{equation}

%24
\begin{equation}
\begin{array}{rcl}
\tR(\vX,\hY)\vZ=\{ [K_Y,K_X]Z -\r(Y,X,Z)-(\nabla K)(Y,X,Z)\}^h,
\end{array}
\end{equation}

%25
\begin{equation}
\tR(\hX,\hY)\hZ=(\gR(X,Y)Z)^h ,
\end{equation}
%26
\begin{equation}
\tR(\hX,\hY)\vZ=(\gR(X,Y)Z)^v.
\end{equation}

\end{corollary}

\begin{corollary}\label{tR_Hessian}
For a Hessian structure the curvature tensor $\tR$ is given by the
formulas
\begin{equation}\label{formulas_tR1}
\begin{array}{rcl}
&&\tR(X^v,Y^v)Z^v=\tR(X^h,Y^h)Z^v=([K_Y,K_X]Z)^v=(\gR(X,Y)Z)^v,\\
&&\tR(X^v,Y^v)Z^h=\tR(X^h,Y^h)Z^h=([K_Y,K_X]Z)^h=(\gR(X,Y)Z)^h,\\
&&\tR(X^v,Y^h)Z^v=-\tR(Y^h,X^v)Z^v\\
&&\ \ \ \ \ \ \ \ \ \ \ \ =-\nabla K(Y,X,Z)^h+([K_Y,K_X]Z)^h\\
&&\ \ \ \ \ \ \ \ \ \ \ \  =-\nabla K(Y,X,Z)^h+(\gR(X,Y)Z)^h,\\
&&\tR(X^v,Y^h)Z^h=-\tR(Y^h,X^v)Z^h\\
&&\ \ \ \ \ \ \ \ \ \ \ \ \ =\nabla K(Y,X,Z)^v-([K_Y,K_X]Z)^v\\
&&\ \ \ \ \ \ \ \ \ \ \ \ \   =\nabla K(Y,X,Z)^v-(\gR(X,Y)Z)^v.
\end{array}
\end{equation}
\end{corollary}

Our task is now to compute the Ricci tensor $\tRic$ of $\tnabla$. In
what follows $e_1,..., e_n$ will stand for an orthonormal basis of
$T_p\M$. As usual $\xi\in T_p\M$ is a fixed point of the tangent
bundle. We have

\begin{equation}\label{for_Ric_hh}
\begin{array}{rcl}
&&\sum_{i=1}^ng(R(Y,\r(Z,\u,e_i))\u,e_i)\\
&&\\
&&\ \ =\sum_{ij=1}^n g(R(Y,g(\r(Z,\u,e_i),e_j)e_j)\u,e_i)\\
&&\\
&&\ \ =\frac{1}{2}\sum_{ij}g(R(e_j,Z)\u,e_i)g(R(Y,e_j)\u,e_i)\\
&&\\
&&\ \ =-\frac{1}{2}\sum_{ij=1}^ng(R(e_j,Z)\u,e_i)g(R(e_j,Y)\u,e_i)\\
&&\\
 &&\ \ =-\frac{1}{2}g(R_2(Z,\xi), R_2(Y,\xi)),
\end{array}
\end{equation}
where the $(1,1)$-tensor $R_2(U,W)$ for two fixed vectors $U,W\in
T_p\M$ is defined by
\begin{equation}\label{R_1}
R_2(U,W)(V)=R(V,U)W.
\end{equation}
For later use we shall also introduce the $(0,2)$-tensor $R_3(V,W)$,
where $V,W\in T_p\M$ given by
\begin{equation}\label{R_2}
R_3(V,W)(X,Y)=g(R(X,Y)V,W)
\end{equation}
and the $(1,2)$-tensor $R_4(V)$, where $V\in T_p\M$, given by
\begin{equation}\label{R_3}
R_4(V)(X,Y)=R(X,Y)V.
\end{equation}

Using  Proposition \ref{curvature_tR}, (\ref{for_Ric_hh}),
(\ref{about_traces}) and (\ref{dtau}) we can now compute
$\tRic(\hY,\hZ)$ at $\u\in T_p\M$. We have

\begin{eqnarray*}
&&\tRic(\hY,\hZ)=\sum_{i=1}^n
\tg(\tR(e_i^v,\hY)\hZ,e_i^v)+\sum_{i=1}^n
\tg(\tR(e_i^h,\hY)\hZ,e_i^h)\\
&&\\
&&\ \ =\RicK(Y,Z)-\frac{1}{2}\tr R(Y,Z)+\sum_{i=1}^ng(\nabla K(Y,Z,e_i),e_i)\\
&&\\
 &&\ \ \ \ -\frac{1}{2}\sum_{i=1}^ng(R(Y,\r(Z,\u,e_i)\u,
e_i))+\gRic
(Y,Z)+\frac{1}{2}\sum_{i=1}^ng(\r(e_i,\u,R(Y,Z)\u),e_i)\\
&&\\
 &&\ \ \ \ \
-\frac{1}{2}\sum_{i=1}^ng(\r(Y,\u,R(e_i,Z)\u)e_i)-\sum_{i=1}^ng(\r(Z,\u,R(e_i,Y)\u),e_i)\\
&&\\
 &&=\RicK (Y,Z)+\gRic(Y,Z)-\frac{1}{2}\tr
R(Y,Z)+\nabla\tau(Y,Z)\\
&&+\frac{1}{4}\sum_{ij=1}^ng(R(e_j,Z)\u,e_i)g(R(e_j,Y)\u,e_i)-\frac{3}{4}\sum_{i=1}^ng(R(e_i,Y)\u,R(e_i,Z)\u)\\
&&\\
 &&=\RicK (Y,Z)+\gRic(Y,Z)-\frac{1}{2}\tr
R(Y,Z)+\nabla\tau(Y,Z)-\frac{1}{2}\sum_{i=1}^ng(R(e_i,Y)\u,R(e_i,Z)\u)\\
&&\\
&&=\RicK (Y,Z)+\gRic(Y,Z)
+\frac{1}{2}\nabla\tau(Y,Z)+\frac{1}{2}\nabla\tau(Z,Y)-\frac{1}{2}g(R_2(Y,\xi),R_2(Z,\xi))\\
&&\\
&& =\frak{Ric} (Y,Z)
+\frac{1}{2}\nabla\tau(Y,Z)+\frac{1}{2}\nabla\tau(Z,Y)-\frac{1}{2}g(R_2(Y,\xi),R_2(Z,\xi)).
\end{eqnarray*}

We shall now compute $\tRic(\vY,\hZ)$. By Proposition
\ref{curvature_tR} we obtain

\begin{equation}\label{ricci_vh_1}
\begin{array}{rcl}
&&\tRic (\vY,\hZ)\\
&&\ \ \ \ =\sum_{i=1}^n\tg(\tR(e_i^v,\vY)\hZ,e_i^v)+\sum_{i=1}^n\tg(\tR(e_i^h,\vY)\hZ,e_i^h)\\
&&\ \ \ \
=\sum_{i=1}^n g(K(\r(Z,\xi,Y),e_i),e_i)-\sum_{i=1}^n g(K(\r(Z,\xi,e_i),Y),e_i)\\
&&\ \ \ \ \ \
+\frac{1}{2}\sum_{i=1}^ng(K(Y,R(e_i,Z)\xi),e_i)-\sum_{i=1}^n
g(\r(K(e_i,Z),\xi,Y),e_i)\\
&&\ \ \ \ \ \  -\sum_{i=1}^n g((\nabla_{e_i} \r) (Z,\xi,Y),e_i)
+\sum_{i=1}^ng(K(e_i,\r(Z,\xi,Y),e_i)\\
&&\ \ \ \ \ \ \ \ \ \ \ \ \ \ \ \ \ \ \ \ \ \ \ \ \ \ \ \ \ \ \ \
 +\sum_{i=1}^ng(\r(e_i,\xi,
K(Y,Z)),e_i)\\
&&\ \ \ \ \ \   =g(R(E,Z)\xi,Y) -\frac{1}{2}g(R(K(Y,e_i),Z)\xi,e_i)\\
&& \ \ \ \  \ \ \ \ \
+\frac{1}{2}\sum_{i=1}^ng(R(e_i,Z)\xi,K(Y,e_i))-\frac{1}{2}g(R(e_i,K(e_i,Z))\xi,Y)\\
&&\ \ \ \ \ \ \ \ \ \ \ \ \ \ \ \ \ \ \ \ \ \ \ \ \ \ \ \ \ \ \ \ \
 -\sum_{i=1}^ng((\nabla_{e_i} \r) (Z,\xi,Y),e_i)\\
&&\ \ \ \ \ \   =g(R(E,Z)\xi,Y) -\frac{1}{2}g(R(K(Y,e_i),Z)\xi,e_i)\\
&& \ \ \ \  \ \ \ \ \
+\frac{1}{2}\sum_{i=1}^ng((K_Y\circ R_2(Z,\xi))e_i,e_i)-\frac{1}{2}g(R(e_i,K(e_i,Z))\xi,Y)\\
&&\ \ \ \ \ \ \ \ \ \ \ \ \ \ \ \ \ \ \ \ \ \ \ \ \ \ \ \ \ \ \ \ \
 -\sum_{i=1}^ng((\nabla_{e_i} \r) (Z,\xi,Y),e_i).
\end{array}
\end{equation}
Since $\tr(K_Y\circ R_2(Z,\xi))=\tr(R_2(Z,\xi)\circ K_Y)$, we have
\begin{equation}\label{ricci_vh_2}
\sum_{i=1}^ng(R(e_i,Z)\xi,K(Y,e_i))=\sum_{i=1}^ng(R(K(Y,e_i),Z)\xi,e_i).
\end{equation}
We also have
\begin{equation}\label{ricci_vh_3}
\sum_{i=1}^nR(e_i,K(e_i,Z))=0.
\end{equation}
Indeed, the following equalities hold
\begin{eqnarray*}
&&\sum_{i=1}^n
R(e_i,K(e_i,Z))=\sum_{i,j=1}^nR(e_i,g(K(e_i,Z),e_j)e_j)\\
&&\ \
=\sum_{i,j=1}^nR(g(K(e_j,Z),e_i)e_i,e_j)=\sum_{j=1}^nR(K(e_j,Z),e_j)=-\sum_{i=1}^nR(e_i,K(e_i,Z)).
\end{eqnarray*}
To complete the computations of $\tRic(Y^v,Z^h)$ we  prove

\begin{lemma}\label{nablaXr}
For any statistical structure we have
\begin{equation}
\begin{array}{rcl}
&&g((\nabla_X\r)(Z,V,Y),W)=\frac{1}{2}g((\nabla_XR)(W,Z)V,Y)\\
&&\ \ \ \ \ \ \ \ \ \ \ +g(R(K(X,W),Z)V,Y)-g(R(W,Z)V,K(X,Y)).
\end{array}
\end{equation}
\end{lemma}
\proof Let $X,Y,Z,V,W\in T_p\M$. Extend the vectors $Y,Z,V,W$ to
vector fields $Y,Z,V,W$ defined in a neighborhood of $p$ in such a
way that $\gnabla_X Y=\gnabla_X Z=\gnabla_X V=\gnabla_X W=0$ at $p$.
We now compute at $p$

\begin{eqnarray*}
&&g((\gnabla_X\r)(Z,V,Y),W)=g(\gnabla_X(\r(Z,V,Y)),W)\\
&&\ \ \ \ \ \ \ =X(g(\r(Z,V,Y),W))=\frac{1}{2}X(g(R(W,Z)V,Y))\\
&&
=\frac{1}{2}g((\gnabla_XR)(W,Z)V,Y)=\frac{1}{2}g((\nabla_XR)(W,Z)V,Y)-\frac{1}{2}g((K_XR)(W,Z)V,Y)\\
&&=\frac{1}{2}g((\nabla_XR)(W,Z)V,Y)-\frac{1}{2}\{g(K(X,Y),
R(W,Z)V)\\
&&\ \ \ \ \ \ \ \
 -g(R(K(X,W),Z)V,Y)-g(R(W,K(X,Z))V,Y)\\
 &&\ \ \ \ \ \ \ \ \ \ \ \ \ \ \ \ \  \ \ \ \ \ \ \ \ \ \ \ \ \ \ -g(R(W,Z)K(X,V),Y)\},
\end{eqnarray*}
\begin{eqnarray*}
&&g((K_X\r)(Z,V,Y),W)=g(\r(Z,V,Y),
K(X,W))-g(\r(K(X,Z),V,Y),W)\\
&&\ \ \ \ \ \ \ \ \ -g(\r(Z,K(X,V),Y),W)
 -g(\r(Z,V, K(X,Y)),W)\\
&&=\frac{1}{2}\{g(R(K(X,W),Z)V,Y)-g(R(W,K(X,Z))V,Y)\\
&&\ \ \ \ \ \  \ \ \  -g(R(W,Z)K(X,V),Y)-g(R(W,Z)V,K(X,Y))\}.
\end{eqnarray*}
Since $\nabla_X\r=\gnabla_X\r+K_X\r$, we obtain the desired formula.
\koniec

Lemma \ref{nablaXr} yields
\begin{equation}\label{ricci_vh_4}
\begin{array}{rcl}
&&\sum_{i=1}^ng((\nabla_{e_i}\r)(Z,\xi,Y),e_i)=\frac{1}{2}\sum_{i=1}^ng((\nabla_{e_i}R)(e_i,Z)\xi,Y)\\
&& \ \ \ \ \ \ \ \ \ \ \ \ \ \ \  \  \  \
+g(R(E,Z)\xi,Y)-g(R_2(Z,\xi),K_Y).
\end{array}
\end{equation}
Combining (\ref{ricci_vh_1})-(\ref{ricci_vh_4}) one gets

\begin{equation}
\tRic(Y^v,Z^h)_{|\xi}=-\frac{1}{2}\tr_gg((\nabla_{\cdot}R)(\cdot,Z)\xi,Y)+g(R_2(Z,\xi),K_Y).
\end{equation}

 In a similar way, using Proposition \ref{curvature_tR}
we can find a formula for $\tRic(\vY,\vZ)$. We have

\begin{eqnarray*}
&&\tRic(\vY,\vZ)=\sum_{i=1}^n\tg(\tR(e_i^v,\vY)\vZ,e_i^v)+\sum_{i=1}^n\tg(\tR(e_i^h,\vY)\vZ,e_i^h)\\
&&\ \ \ \ =-2\RicK(Y,Z)+\sum_{i=1}^ng(\r(e_i,Y,Z),e_i)\\
&&\ \ \ \ \ \ \ +\sum_{i=1}^ng(\nabla
K(e_i,Y,Z),e_i)-\sum_{i=1}^ng(\r(\r(e_i,\u,Z),\u,Y),e_i)\\
&&\ \ \ \ =-2\RicK(Y,Z)+\sum_{i=1}^ng(\nabla
K(e_i,Y,Z),e_i)-\frac{1}{2}\sum_{ij=1}^ng(R(e_i,g(\r(e_i,\u,Z),e_j)e_j)\u,Y)\\
&&\ \ \ \ =-2\RicK(Y,Z)+\sum_{i=1}^ng(\nabla
K(e_i,Y,Z),e_i)-\frac{1}{4}\sum_{ij=1}^ng(R(e_i,e_j)\u,Y)g(R(e_j,e_i)\u,Z)\\
&&\ \ \ \
=\frac{1}{2}(\Ric(Y,Z)-\Ricstar(Y,Z))+\nabla\tau(Y,Z)-\frac{1}{4}g(R_3(\xi,Y),R_3(\xi,Z)).
\end{eqnarray*}
For the last equality we have used (\ref{for_Ric_vv}) and
(\ref{Ric+Ricstar}).

Summing up,  we have proved

\begin{proposition}\label{Ricci_tg}
For any statistical structure and the induced Riemannian  structure
on the tangent bundle we have the following formulas for the Ricci
tensor $\tRic$  at $\xi\in T_p\M$ for every $Y,Z\in T_p\M$
\begin{equation}\label{tRic_vv}
\begin{array}{rcl}
&&\tRic (\vY,\vZ)= -2\RicK(Y,Z)+(\div ^\nabla K)(Y,Z)  \\
&& \ \ \ \ \ \ \ \ \ \ \ \ \ \ \ \ \ \ \ \ \ \ \ \ \ \
-\frac{1}{4}g(R_3(\xi,Y),R_3(\xi,Z))\\
&&\ \ \ \ \ \ \ \ \ \ \ \ =\frac{1}{2}(\Ric(Y,Z)-\Ricstar
(Y,Z))+\nabla \tau
(Y,Z)\\
&&\ \ \ \ \ \ \ \ \ \ \ \ \ \ \ \ \ \ \ \ \ \ \ \ \ \ \
 -\frac{1}{4}g(R_3(\xi,Y),R_3(\xi,Z)),\\
 %&&\ \ \ \ \ \ \ \ \ \ \ \ \ \ \ \ \  =-(\div^{\gnabla}K)(Y,Z)-\tau(K(Y,Z))\\
%&& \ \ \ \ \ \ \ \ \ \ \ \ \ \ \ \ \ \ \ \ \ \ \ \ \ \ \
% -\frac{1}{4}g(R_3(\xi,Y),R_3(\xi,Z)),
\end{array}
\end{equation}
\begin{equation}\label{tRic_vh}
\gRic(Y^v,Z^h)=-\frac{1}{2}\tr_gg((\nabla_{\cdot}R)(\cdot,Z)\xi,Y)+g(R_2(Z,\xi),K_Y),
\end{equation}
\begin{equation}\label{tRic_hh}
\begin{array}{rcl}
&&\tRic(\hY,\hZ)=\frak{Ric} (Y,Z))+\frac{1}{2}\nabla\tau(Y,Z)+\frac{1}{2}\nabla\tau(Z,Y)\\
&&\ \ \ \ \ \ \ \ \ \ \ \ \ \ \ \ \ \ \ \ \ \ \ \ \ \ \ \ \
 -\frac{1}{2}g(R_2(Y,\xi),R_2(Z,\xi)),
\end{array}
\end{equation}
where $R_2$, $R_3$ are defined by {\rm(\ref{R_1})} and
{\rm(\ref{R_2})}.
\end{proposition}

\begin{proposition}\label{prop_scalar_for_tg}
For any statistical structure we have the following formula for the
scalar curvature $\rhotg_\u$ of the  Sasakian  metric on the tangent
bundle
\begin{equation}\label{scalar_for_tg}
\rhotg_\u=\rho+2\tr_g\nabla \tau-\frac{3}{4}\Vert R_4(\xi)\Vert^2,
%\sum_{ij=1}^n\VertR(e_i,e_j)\u\Vert^2.
\end{equation}
where $R_4$ is defined by {\rm(\ref{R_3})}.
\end{proposition}

Denote by  $\tk(X^a\wedge Y^b)$, where $a,b=v $ or $h$, the
sectional curvature of $\tg$ by a plane spanned by $X^a,Y^b$.

Let $X^a,Y^b$ be a pair of unit vectors in $T_{\u}T\M$. By
Proposition \ref{curvature_tR} we have
\begin{proposition}\label{prop_sectional_curvature}
For any statistical structure we have
\begin{equation}\label{sectional_vv}
\tk(\vX\wedge\vY)=-\Kk (X\wedge Y),
\end{equation}
\begin{equation}\label{sectional_hh}
\tk(\hX\wedge\hY)=\gk (X\wedge Y)-\frac{3}{4}\Vert R(X,Y)\u\Vert^2,
\end{equation}
\begin{equation}\label{sectional_vh}
\tk (\vX\wedge\hY)=g([K_X,K_Y]Y,X)+g(\nabla K(Y,X,Y),X)+\Vert
\r(Y,\u,X)\Vert ^2.
\end{equation}
If the statistical structure is conjugate symmetric, then the last
formula becomes
\begin{equation}
\tk (\vX\wedge\hY)=g([K_X,K_Y]Y,X)+g(\nabla
K(Y,X,Y),X)+\frac{1}{4}\Vert R(X,\u)Y\Vert ^2
\end{equation}

\end{proposition}

\proof The first two equalities immediately  follow from Proposition
\ref{curvature_tR}.  We shall check (\ref{sectional_vh}). First
observe that
\begin{eqnarray*}
&&g(R(Y,\r(Y,\u,X))\u,X)=\sum_{i=1}^ng(R(Y,g(\r(Y,\u,X),e_i)e_i)\u,X)\\
&&\ \ \ \ \ =-\frac{1}{2}\sum_{i=1}^ng(R(Y,e_i)\u,X)^2=-2\Vert
\r(Y,\u,X)\Vert^2.
\end{eqnarray*}
Thus
\begin{eqnarray*}
&&\tg(\tR(\vX,\hY)\hY,\vX)=g([K_X,K_Y]Y,X)+g(\nabla
K(Y,X,Y),X)\\
&&\ \ \ \ \ \ \ \ \ \ \ \ \ \ \ \ \ \ \ \ \ \ \ \ \ \ \
 -\frac{1}{2}g(R(Y,\r(Z,\u,X))\u,X)\\
&&\ \ \ \ \  =g([K_X,K_Y]Y,X)+g(\nabla K(Y,X,Y),X)+\Vert
\r(Y,\u,X)\Vert^2.
\end{eqnarray*}
\koniec

Recall the classical theorem of Sasaki, \cite{Sasaki}.
\begin{thm}
Let $g$ be a metric tensor on $\M$ and $\tg$ be the Sasaki metric on
$T\M$ determined by $g$ and $\gnabla$. Then $\tg$ is flat if and
only if $g$ is flat.
\end{thm}
In \cite{Gud} a few results generalizig this theorem were proved.
For instance,

\begin{thm} Let $g$ be a metric tensor on $\M$ and $\tg$  the Sasaki metric on
$T\M$ determined by $g$ and $\gnabla$. If the sectional  curvature
of $\tg$ is  bounded on $T\M$, then $g$ is flat on $\M$.
\end{thm}

\begin{thm}\label{thm_scalar_curv_tan}
Let $g$ be a metric tensor on $\M$ and $\tg$ be the Sasaki metric on
$T\M$ determined by $g$ and $\gnabla$. The scalar curvature of $\tg$
is constant on $T\M$ if and only if $g$ is flat on $\M$.
\end{thm}

 By Propositions \ref{prop_scalar_for_tg}  \ref{prop_sectional_curvature} one  gets the following
 generalizations of the above
 theorems to the case od statistical structures.

\begin{thm}\label{thm_sec_curv_tan} Let $(g,\nabla)$ be a statistical structure on $\M$ and
$\tg$ be the Sasakian metric on $T\M$ determined by $g$ and
$\nabla$. If $\nabla$ is not flat, then the sectional curvature of
$\tg$ on $T\M$ is unbounded from above and below.
\newline
If $\nabla$ is not flat, then the scalar curvature of $\tg$ is
unbounded from below on $T\M$.
\end{thm}

Note that in the case of the Sasakian metrics determined by
statistical structures, the flatness of a statistical connection
does not imply the flatness of the Sasakian metric. Namely, we have

\begin{thm}\label{thm_flatness_tan}If the Sasakian metric on the tangent bundle of a statistical manifold is
flat, then the statistical structure is Hessian, its Riemannian
metric is flat and its difference tensor is parallel relative to the
statistical connection. The converse is also true.
\end{thm}
\proof  Assume that the Sasakian metric is flat. From the first two
formulas in Proposition \ref{prop_sectional_curvature} we see that
the curvature tensors $[K,K]$,  $R^g$ and $R$ vanish. The structure
is, in particular, Hessian. It is now sufficient to look at
Corollary \ref{tR_Hessian}. The converse follows immediately from
Corollary \ref{tR_Hessian}.\koniec

The structure  from the last theorem, although it satisfies a few
strong conditions, does not have to be trivial, that is, it does not
have to reduce to the flat Riemannian structure. Namely, consider
the following example

\begin{example}{\rm Let $\M=\{ x=(x_1,...,x_n)\in \R^{n}:\ \ \ \
x_i>0\ \forall i=1,...,n\}$. Let $g$ be the standard flat metric
tensor field on $\M$ and $e_1,...,e_n$ be the canonical frame field
on $\M$. Define  a symmetric $(1,2)$-tensor field $K$ on $\M$ as
follows

\begin{equation}
\begin{array}{rcl}
 &&K(e_i,e_j)=0 \ \ \  \ for \ \ \ i\ne j,\\
&&K(e_i,e_i)=\lambda_ie_i \ \ \ \  for \ \ \  i=1,...,n,
\end{array}
\end{equation}
where $\lambda_i(x)=-x_i^{-1}$. One sees that $K$ is symmetric
relative to $g$ and $[K,K]=0$. When we define $\nabla=\gnabla+K$,
then the statistical structure $(g,\nabla)$ is non-trivial, Hessian,
$g$ is flat and $\nabla K=0$.
 }
\end{example}

\bigskip

\section{The sphere bundles over statistical manifolds}
\subsection{The sphere bundle as a hypersurface}
%\subsection{The sphere bundle as a hypersurface}
 Denote by $\Sr$ the tangent sphere bundle of radius $r>0$
over $\M$
 equipped with the metric $\tg$ induced by the Sasaki metric
$\tg$ on $T\M$ (determined by a statistical structure on $\M$). Note
that $\Sr$ is a hypersurface in $T\M$. It is clear that the
curvatures of $\Sr$ will depend on the radius $r$. The first task is
to find the normal vector and the tangent space to $\Sr$ at $\xi\in
\Sr$.

Let $\Uv$ denote the canonical vertical vector field on $T\M$, that
is,
$$\Uv_{\xi}=\xi^v_{\xi}.$$
This vector field is not a vertical lift of any vector field on
$\M$. But if $(u^i)$ is a local coordinate system on $\M$, then
\begin{equation}\label{U_in_coord}
\Uv=(u^i\circ\pi)\partial_i^v.
\end{equation}
In the following lemma some basic information about the canonical
vertical  vector field is collected.
\begin{lemma}\label{diff_Uv}
  For any vector fields $X,Y$ on $\M$ we have
\begin{equation}
(\tnabla_{X^v}\Uv)_\xi=K(X,\xi)^h_\xi+X^v_\xi,
\end{equation}
\begin{equation}\label{XhUv}
(\tnabla_{X^h}\Uv)_\xi=-\r(X,\xi,\xi)^h_\xi-K(X,\xi)^v_\xi,
\end{equation}
\begin{equation}
(\tnabla_{\Uv} X^v)_\xi=K(X,\xi)^h_\xi,
\end{equation}
\begin{equation}
(\tnabla_{\Uv} X^h)_\xi=-\r(X,\xi,\xi)^h_\xi-K(X,\xi)^v_\xi,
\end{equation}
\begin{equation}
(\tnabla_{\Uv}\Uv)_\xi=K(\xi,\xi)^h_\xi+\U_\xi.
\end{equation}
\end{lemma}
\proof The first two formulas immediately follow from Lemma
\ref{tensor_of_type11} for $\F=\id$. As concerns the two subsequent
formulas we have, see (\ref{U_in_coord}),
(\ref{lifts_with_function})
 and Lemma \ref{tnabla},

\begin{equation*}
(\tnabla_{\Uv}
X^v)_\xi=(u^i\circ\pi)\tnabla_{\partial^v_i}X^v=(u^i\circ\pi)K(\partial_i,X)^h_\xi=K(\xi,X)^h_\xi,
\end{equation*}
\begin{equation*}
\begin{array}{rcl}
&&(\tnabla_{\Uv}
X^h)_\xi=(u^i\circ\pi)\tnabla_{\partial^v_i}X^h=(u^i\circ\pi)[-K(\partial_i,X)^v_\xi-\r(X,\xi,\partial_i)^h_\xi]\\
&&\ \ \ \ \ \ \ \ \ \ \ \ \ \ \ \ \ \ \ \
=-K(\xi,X)^v_\xi-\r(X,\xi,\xi)^h_\xi.
\end{array}
\end{equation*}
To check the last formula we write

\begin{equation*}
\begin{array}{rcl}
&&(\tnabla_{\Uv}\Uv)_\xi=((u^i\circ\pi)\tnabla_{\partial_i^v}\Uv)_\xi\\
&&\ \ \ \ \ \ \ \ \ \ \ =(u^i\circ
\pi)K(\partial_i,\xi)^h_\xi+(u^i\circ\pi)(\partial_i^v)_\xi=K(\xi,\xi)^h_\xi+\xi^v_\xi.
\end{array}
\end{equation*}\koniec

We shall now find vectors spanning the tangent space and the normal
vector at $\xi\in\Sr$.
\begin{lemma}\label{lemma_tangentspace_N}
If $X\in T_p\M$ and $\xi\in \Sr$, then the vectors

\begin{equation}
X^{tv}_\xi:=X^v_{\xi}-\frac{1}{r^2}g(\xi,X)\U_\xi
\end{equation}
and
\begin{equation}
X^{th}_\xi:=X^h_\xi +\frac{1}{r^2}g(K(\xi,\xi),X)\U_\xi
\end{equation}
are tangential to $\Sr$ at $\xi$.

The normal vector field $\N$ at $\xi\in T_p\M$ is given by
\begin{equation}\label{N}
\N_\xi=\frac{1}{\sqrt{r^2+\Vert K(\xi,\xi)\Vert^2}}(\Uv_\xi
-(K(\xi,\xi))^h).
\end{equation}
\end{lemma}
\proof Let $X\in T_p\M$ and $\xi_t$ be the integral curve (as in the
previous sections) of a  lift (vertical or horizontal) of $X$ to
$\xi\in \Sr$. Set $\psi(t):=\Vert \xi_t\Vert$ and
$\phi(t):=\psi(t)^{-1}$. The vectors $r[t\to \phi(t)\xi_t]$ are
tangent to $\Sr$ and they span $T_\xi(\Sr)$.

Consider first the vertical lift, that is, $\xi_t=\xi+tX$. Then
$\psi(t)^2=r^2+t^2\Vert X\Vert^2+2tg(\xi,X)$ and thus $\vX
(\psi(t)^2)=2g(\xi,X)$. Thus $\phi'(0)=-r^{-3}g(\xi,X)$. We have
$r[t\to \phi(t)\xi_t]=X^v_{\xi}-\frac{1}{r^2}g(\xi,X)\xi^v_\xi$.
Recall that if $\xi$ is treated as a vector, then it can be
identified with $\xi^v$. Let now $\xi_t$ be the integral curve of
$\hX$ defined in the previous part of this paper, that is,
$\nabla_X\xi_t=0$. By (\ref{XhUv})
 we have

\begin{equation*}
\hX(\psi(t))^2=\hX(\tg(\Uv,\Uv))=2\tg(\tnabla_{\hX}\Uv,\Uv)=-2g(K(\xi,\xi),X).
\end{equation*}
 It follows that $\phi'(0)=r^{-3}g(K(\xi,\xi),X)$, and consequently,
$r[t\to \phi(t)\xi_t]=X^h_\xi
+\frac{1}{r^2}g(K(\xi,\xi),X)\xi^v_\xi$.

One now easily checks that $\N$ given by (\ref{N}) is orthogonal to
the vectors $X^{tv}$ and $X^{th}$.
 \koniec

 The vectors $X^{tv}$, $X^{th}$ will be called  tangent
vertical and
 tangent horizontal respectively. The assignments $\T_p\M \ni X\to X^{tv}$,
 $X\to X^{th}$ are linear.

The tangent space $T_\xi(\Sr)$  can be regarded  as the direct sum
of the two subspaces

\begin{equation}\label{Vr_Hr}
{\mathcal V}_\xi^r:=\{X^v,\  where \ X\bot\xi \}\oplus {\mathcal
H}_\xi ^r:=\{X^\th \ \ for\ X\in T_p\M\}
\end{equation}
or the direct sum of the three subspaces
\begin{equation}
\{X^v,\ where \ X\bot\xi\}\oplus \{X^h,\  where \ X\bot
K(\xi,\xi)\}\oplus span\{K(\xi,\xi)^{th}\}.
\end{equation}
The union  $\mathcal V^r=\bigcup _{\xi\in\Sr}\mathcal V^r_\xi$ forms
a smooth distribution over $\Sr$.
 Note that the mapping $T_p\M\ni X\to X^{th}\in\mathcal H_\xi^r$ is not an
 isometry, in general. Of course, the mapping $\{\xi\}^\bot\ni X\to
 X^{th}$ is an isometry.

In what follows $\h$ will denote  the second fundamental form for
the hypersurface $\Sr$ in $T\M$.

\begin{proposition}\label{h}

Let $\xi \in\Srp$ and $X,Y\in T_p\M$.
\newline
For every  $X,Y$ we have

\begin{equation}\label{h_thth}
\begin{array}{rcl}
&&\h(X^\th_\xi,Y^\th_\xi)\\
&& \\
&&=\f\{\frac{1}{2}[g((\nabla^g_XK)(\xi,\xi),Y)+g((\nabla^g_YK)(\xi,\xi),X)]\\
&& \\
 &&\ \ \ \ \ \ \ \
-\frac{1}{r^2}[g(K(\xi,\xi),Y)g(R(X,K(\xi,\xi))\xi,\xi)\\
&&\ \ \ \ \ \ \ \ \ \ \ \ \ \ +g(K(\xi,\xi),X)g(R(Y,K(\xi,\xi))\xi,\xi)]\\
&&\ \ \ \ \ \ \ \ \ -2g(K(\xi,X),K(\xi,Y))\\
&& \\
&&\ \ \ \ \ \ \ \ \
+\frac{1}{r^2}g(K(\xi,\xi),X)g(K(\xi,\xi),Y)(1-\frac{1}{r^2}\Vert
K(\xi,\xi)\Vert^2)\}.
\end{array}
\end{equation}
For every $X$ and $Y\bot \, \xi$ we have

\begin{equation}\label{h_vth}
\begin{array}{rcl}
&&\h(X^\th_\xi,Y^v_\xi)\\
&& \\
&&\ \ \ \ =\f[g(K(X,Y),\xi)-\frac{1}{2}g(R(X,K(\xi,\xi))\xi,Y)\\
&& \\
 &&\ \ \ \ \ \ \ \
-\frac{1}{r^2}g(K(\xi,\xi),X)g(K(\xi,Y),K(\xi,\xi))].
\end{array}
\end{equation}
For every $X\bot\xi$ and every  $Y$ we have

\begin{equation}\label{h_vth}
\begin{array}{rcl}
&&\h(\vX_\xi,Y^{th}_\xi)\\
&& \\
&&\ \ \ \ =\f[g(K(X,Y),\xi)-\frac{1}{2}g(R(Y,K(\xi,\xi))\xi,X)\\
&& \\
 &&\ \ \ \ \ \ \ \
-\frac{1}{r^2}g(K(\xi,\xi),Y)g(K(\xi,X),K(\xi,\xi))].
\end{array}
\end{equation}
 For every $X,Y\bot\,\xi$ we have

\begin{equation}\label{h_vv}
\h(X^v_\xi,Y^v_\xi)=\f[-g(X,Y)-g(K(X,Y),K(\xi,\xi))].
\end{equation}
\end{proposition}

\proof In the proof we shall use the notations:
$f(\u)=\frac{1}{\sqrt{r^2+K(\u,\u)}}$, $\N_1(\u)=\Uv
_{\u}-K(\u,\u)^h_\xi$. Thus $\N=f\N_1$ and clearly
$\tg(\N_1,Y^{th})=0=\tg(\N_1,Y^{tv})$. The identity mapping on $T\M$
will be denoted by $u$.

By (\ref{tensor_K_vh}) and (\ref{tensor_K_hh}) we have
\begin{equation}\label{tnabla_hKuuh}
\tnabla_{X^h_\xi}(K(u,u))^h=-\frac{1}{2}(R(X,K(\xi,\xi))\xi)^v_\xi+((\gnabla_XK)(\xi,\xi))^h_\xi-2(K^2_\xi(X))^h_\xi,
\end{equation}

\begin{equation}\label{vKuu_h}
\tnabla_{X^v_\xi}
(K(u,u))^h=-K(X,K(\xi,\xi))^v_\xi-\r(K(\xi,\xi),\xi,X)^h_\xi
+2K(X,\xi)^h_\xi.
\end{equation}
In particular,

\begin{equation}\label{tnabla_Kuuh}
\tnabla_{\xi^v_\xi}
(K(u,u))^h=-K(\xi,K(\xi,\xi))^v_\xi-\r(K(\xi,\xi),\xi,\xi)^h_\xi
+2K(\xi,\xi)^h_\xi.
\end{equation}

In the following computations made at $\xi$ we shall skip $\xi$ in
the subscript position. Using Lemma \ref{diff_Uv} and
(\ref{tnabla_hKuuh}) we obtain
\begin{equation*}
\begin{array}{rcl}
&&\tnabla_{X^h}(f\N_1)=(X^hf)\N_1+f(\xi)[-\r(X,\xi,\xi)^h-K(X,\xi)^v ]\\
&&\ \ \ \ \ \ \ \ \ \
+f(\xi)[\frac{1}{2}(R(X,K(\xi,\xi)\xi)^v)-((\gnabla_XK)(\xi,\xi))^h+2(K^2_\xi
X)^h)].
\end{array}
\end{equation*}
 Therefore, using also (\ref{for_sec_h}), we get
 \begin{equation*}
 \begin{array}{rcl}
&&g(\tnabla_{X^h}\N,Y^{th})=f(\xi)\{-\frac{1}{2}[\{g((\gnabla_XK)(\xi,\xi),Y)+
g((\gnabla_YK)(\xi,\xi),X)]\\
&&\ \ \ \ \ \ \ \ \ \ \ \ \ \ \ \ \ \ \ \ \ \ \ \ \ \ \ \ \ \ \ \ -\frac{1}{r^2}g(K(\xi,\xi),X)g(K(\xi,\xi),Y)\\
&&\ \ \ \ \ \ \
+\frac{1}{2r^2}g(R(X,K(\xi,\xi))\xi,\xi)g(K(\xi,\xi),Y)+2g(K(\xi,X),K(\xi,Y))\}.
\end{array}
 \end{equation*}
Similarly, making straightforward computations (using Lemma
\ref{diff_Uv} and (\ref{tnabla_Kuuh})), we receive

\begin{equation*}
\begin{array}{rcl}
&&\tnabla_{\xi^v}(f\N_1)=(\xi^vf)\N_1+f(\xi)[K(\xi,\xi)^h+\xi^v]\\
&&\ \ \ \ \ \ \ \ \ \ \ \ \ \
-f(\xi)[2K(\xi,\xi)^h-\r(K(\xi,\xi),\xi,\xi)^h-K(\xi,K(\xi,\xi))^v],
\end{array}
\end{equation*}
and consequently,
\begin{equation*}
\begin{array}{rcl}
&&g(\tnabla_{\frac{1}{r^2}g(K(\xi,\xi),X)\xi^v}\N, Y^{th})\\
&&\ \ \ \ \ \ \ \ \ \ \ \
=f(\xi)\{\frac{1}{2r^2}g(R(Y,K(\xi,\xi)\xi,\xi)g(X,K(\xi,\xi))\\
&&\ \ \ \ \ \ \ \ \ \ \ \ \ \ \ \ \ \ \ \ \ \ +\frac{1}{r^4}\Vert
K(\xi,\xi)\Vert^2g(K(\xi,\xi),X)g(K(\xi,\xi),Y)\}.
\end{array}
\end{equation*}
We have proved (\ref{h_thth}). In order to prove the two other
formulas we take $X\in T_{\pi (\xi)}\M$ orthogonal to $\xi$. Using
Lemma \ref{diff_Uv} and (\ref{vKuu_h}) we obtain

\begin{equation*}
\begin{array}{rcl}
&&\tnabla_{X^v}(f\N_1)=(X^vf)\N_1+f(\xi)[K(X,\xi)^h+X^v]\\
&&\ \ \ \ \ \ \ -f(\xi)[2K(X,\xi)^h
-K(X,K(\xi,\xi))^v-\r(K(\xi,\xi),\xi,X)^h].
\end{array}
\end{equation*}
One now straightforwardly computes $\tg(\tnabla_{X^v}(f\N_1),Y^v)$
for $Y\bot\xi$ and $\tg(\tnabla_{X^v}(f\N_1),Y^{th})$ for any $Y\in
T_{\pi (\xi)}\M$.\koniec

Having  the above formulas for the second fundamental form we can
compute the mean curvature.

\begin{proposition}\label{H}
 The mean curvature $H$ for $\Sr$ in $TM$ is
given by the following formula
\newline
\begin{equation}\label{mean_curvature}
\begin{array}{rcl}
&&H_{\xi}=\frac{1}{\sqrt{r^2+\Kxixi^2}}\{-(n-1)-\tau(\Kxi)+\frac{\Vert\Kxi\Vert^2}{r^2}\\
&&\\
 &&+\tr(\gnabla_\cdot K)(\xi,\xi)-
\frac{1}{r^2+\Kxixi^2}g((\gnabla_{K(\xi,\xi)}K)(\xi,\xi),\Kxi)\\
&&\\
 &&\ \ \ \ \ \  -2\tr K_\xi^2+2\frac{\Vert
K(\xi,\Kxi)\Vert^2}{r^2+\Kxixi^2}\\
&&\\
 &&\ \ \ \ \ \ \ \ \ \ \ \ \ \
+\frac{1}{r^2+\Kxixi^2}\Kxixi^2(1-\frac{1}{r^2}\Kxixi^2)\}.
\end{array}
\end{equation}

\end{proposition}

\proof We  choose an orhonormal basis of $T_\xi(\Sr)$ as follows. As
usual, let $p=\pi(\xi)$. First we take any orthonormal basis $\tilde
e_1,...,\tilde e_{n-1}$ of the orthogonal complement to $\xi$ in
$T_p\M$ and their vertical lifts $\tilde e_1^v,..., \tilde
e_{n-1}^v$. The vectors $\tilde e_1,...,\tilde
e_{n-1},\frac{\xi}{r}$ form an orthonormal basis of $T_p\M$. If
$\Kxi=0$, then it is sufficient to take any orthonormal basis
$e_1,...,e_n$ of $T_p\M$ and their horizontal lifts
$e_1^h,...,e_n^h$. The set of vectors $\tilde e_1^v,...,\tilde
e_{n-1}^v, e_1^h,..., e_n^h$ is an orhonormal basis of $T_\xi(\Sr)$
in this case. If $\Kxi\ne 0$, then by choosing any orthonormal basis
$e_1,...,e_{n-1}$ of the orthogonal complement to $\Kxi$ in $T_p\M$
we get an orthonormal basis of $T_\xi(\Sr)$ given by $\tilde
e_1^v,...,\tilde e_{n-1}^v, e_1^h,..., e_{n-1}^h, e^*_n,$ where
\begin{equation}
e_n^*=\frac{\Kxi^{th}}{\Vert\Kxi^{th}\Vert}.
\end{equation}
Note that $e_1,...,e_{n-1}, \frac{\Kxi}{\Kxixi}$ is an orthonormal
basis of $T_p\M$.

We shall use the same function  $f$ as in the proof of Proposition
\ref{h}. We have
\begin{equation}
\Kxixit=\frac{\Kxixi}{f(\xi)r}.
\end{equation}
Assume first that $\Kxi\ne 0$.
 We have
 $\h(e^*_n,e^*_n)=\frac{f(\xi)^2r^2}{\Kxixi^2}\h(\Kxit,\Kxit)$,
and consequently (by (\ref{h_thth})),
\begin{equation}
\begin{array}{rcl}
&&\h(e^*_n,e^*_n) =\frac{(f(\xi))^3r^2}{\Kxixi^2}\{g((\gnabla_{\Kxi}K)(\xi,\xi),\Kxi)\\
&&\ \ \ \ \ \ \ \ -2\Vert(K(\xi,\Kxi)\Vert^2
+\frac{1}{r^2}\Kxixi^4(1-\frac{1}{r^2}\Kxixi^2)\}
\end{array}
\end{equation}
and
\begin{equation}
\h(e_i^h,e_i^h)=f(\xi)\{g(\gnabla_{e_i}K)(\xi,\xi),e_i)-2g(K_{\xi}^2e_i,e_i))\}
\end{equation}
for $i=1,...,n-1$. We now obtain

\begin{equation*}
\begin{array}{rcl}
&&\sum_{i=1}^{n-1}\h(e_i^h,e_i^h)+h(e_n^*,e_n^*)\\
&&\\
&&\ \ \ \  \ =f(\xi)\{\tr(\gnabla_\cdot K)(\xi,\xi)\\
&&\\
 &&\ \ \
+\left[\frac{f(\xi)^2r^2}{\Kxixi^2}g\left((\gnabla_{\Kxi}K)(\xi,\xi),\Kxi\right)-
g\left(\gnabla_{\frac{\Kxi}{\Kxixi}}K)(\xi,\xi),\frac{\Kxi}{\Kxixi}\right)\right]\\
&&\\
&& \ \ \ \ -2\tr K_\xi^2+2\left[\Vert
K(\xi,\frac{\Kxi}{\Kxixi})\Vert^2-\frac{f(\xi)^2r^2}{\Kxixi^2}\Vert K(\xi,\Kxi)\Vert^2\right]\\
&&\\
 &&\ \ \ \ \ \ \
+\frac{1}{r^2+\Kxixi^2}\Kxixi^2(1-\frac{1}{r^2}\Kxixi^2)\}\\
&&\\
&&\ \ \ \  \ =f(\xi)\{\tr(\gnabla_\cdot K)(\xi,\xi)-
\frac{1}{r^2+\Kxixi^2}g((\gnabla_{K(\xi,\xi)}K)(\xi,\xi),\Kxi)\\
&&\\
&&\ \ \ \ \ \  +2\tr K_\xi^2-2\frac{\Vert
K(\xi,\Kxi)\Vert^2}{r^2+\Kxixi^2}+\frac{1}{r^2+\Kxixi^2}\Kxixi^2(1-\frac{1}{r^2}\Kxixi^2)\}.
\end{array}
\end{equation*}

If $K(\xi,\xi)=0$, then the same formula is valied, which again
follows from (\ref{h_thth}). Using (\ref{h_vv}) completes the
proof.\koniec

\medskip

\subsection{The sphere bundles with small radii}
Having the explicit formula for the mean curvature of $\Sr$ in $T\M$
we can prove

\begin{thm}\label{thm_H}
Let $(\M, g,\nabla)$  be a compact statistical manifold. For an
arbitrary  number $L$ there exists a sufficiently small radius $r$
such that the absolute value of the mean curvature $H$ of the sphere
bundle $\Sr$ is grater than $L$ everywhere on $\Sr$.
\end{thm}

\proof Assume that the radius $r$ tends to $0$. We take $\lambda r$
for $\lambda \in \R _{+}$ and assume that $\lambda \to 0$. The
initial radius $r$ is fixed. When $r$ is multiplied by $\lambda$,
then, automatically, so is $\xi$. Using (\ref{mean_curvature}) and
making suitable reductions one sees that the function
$$P_1:(0,1]\times \Sr )\ni (\lambda,\xi)\to H_{\lambda\xi}+\frac{n-1}{\sqrt{(\lambda r)^2+\Vert
K(\lambda\xi,\lambda\xi)\Vert^2}}$$  is well defined and continuous
on the compact set $[0,1]\times\Sr$. Hence it is bounded on
$(0,1]\times \Sr$. Since $K(\xi,\xi)$ is bounded on $\Sr$, the
$1$-parameter family of  functions
$$F_\lambda (\xi):(0,1]\times\Sr\ni (\lambda,\xi)\to\frac{n-1}{\lambda\sqrt{
r^2+\lambda^2\Vert K(\xi,\xi)\Vert^2}}$$ tends uniformly to infinity
if $\lambda$ tends to zero, which proves the assertion.\koniec

\begin{thm}\label{thm_normh}
Let $(\M, g,\nabla)$  be a compact statistical manifold. For an
arbitrary  number $L$ there exists a sufficiently small radius $r$
such that the norm of the second fundamental form of the sphere
bundle $\Sr$ in $T\M$ is grater than $L$ everywhere on $\Sr$.
\end{thm}

\proof As in the previous proof we fix  a radius $r$. We use the
same orthonormal bases  as in the proof of Proposition \ref{H}. By
Proposition \ref{h} we get the following expressions for the second
fundamental form $\h$ at $\xi\in \Sr$. If $K(\xi,\xi)\ne 0$, then

\begin{equation}\label{h_nn}
\begin{array}{rcl}
&&\h_\xi(e_n^*,e_n^*)=\frac{r^2}{\Kxixi^2\left(r^2+\Kxixi^2\right)^{\frac{3}{2}}}
\{g((\gnabla_{\Kxi}K)(\xi,\xi),\Kxi)\\
&&\\
 &&\ \ \ \ \ \ \ \ \ \ \ \ \ \ \ -2\Vert(K(\xi,\Kxi)\Vert^2
+\frac{1}{r^2}\Kxixi^4(1-\frac{1}{r^2}\Kxixi^2)\}\\
&&\\
&&\ \ \ \ \ \ \ \
 \ \ \ \ =\frac{r^3}{\left(r^2+\Kxixi^2\right)^{\frac{3}{2}}}
\{g\left((\gnabla_{\frac{\Kxi}{\Kxixi}}K)(\frac{\xi}{r},\xi),\frac{\Kxi}{\Kxixi}\right)\\
&&\\
 &&\ \ \ \ \ \ \ \ \ \ \ \ \ \ \  -2g(K(\frac{\xi}{r},\frac{\Kxi}{\Kxixi}),
 K(\xi,\frac{\Kxi}{\Kxixi}))\\
 &&\\
&& \ \ \ \ \ \ \ \ \ \ \ \ \ \ \
+g(K(\frac{\xi}{r},\xi),K(\frac{\xi}{r},\frac{\xi}{r}))(1-g(K(\xi,\xi),K(\frac{\xi}{r},\frac{\xi}{r}))\},
\end{array}
\end{equation}

\bigskip
\begin{equation}\label{h_in}
\begin{array}{rcl}
&&h_\xi( e_i^h,
e_n^*)=\frac{r}{\Kxixi(r^2+\Kxixi^2)}\{\frac{1}{2}[g((\gnabla_{e_i}K)(\xi,\xi),K(\xi,\xi))\\
&&\\
&&\ \ \ \ \ \ \ \ \ \ \ \ \ \ \ \ \ \ \ \
+g((\gnabla_{K(\xi,\xi)}K)(\xi,\xi),e_i)]\\
&&\\
&&-\frac{1}{r^2}\Kxixi^2g(R(e_i,\Kxi)\xi,\xi)-2g(K(\xi,e_i),
K(\xi,\Kxi))\}\\
&&\\
&&=\frac{r^2}{r^2+\Kxixi^2}\{\frac{1}{2}[g((\gnabla_{e_i}K)(\frac{\xi}{r},\xi),\frac{K(\xi,\xi)}{\Kxixi})\\
&&\\
&&\ \ \ \ \ \ \ \ \ \ \ \ \ \ \ \ \ \ \ \
+g((\gnabla_{\frac{K(\xi,\xi)}{\Kxixi}}K)(\frac{\xi}{r},\xi),e_i)]\\
&&\\
&&-\Kxixi
g(R(e_i,K(\frac{\xi}{r},\xi))\frac{\xi}{r},\frac{\xi}{r})-2g(K(\frac{\xi}{r},e_i),
K(\xi,\frac{\Kxi}{\Kxixi}))\},\\
\end{array}
\end{equation}
\bigskip

\begin{equation}\label{h_tin}
\begin{array}{rcl}
&&\h_\xi(\tilde {e_i}^v,e_n^*)
=\frac{r}{\Kxixi(r^2+\Kxixi^2)}g(K(\tilde e_i,\xi),K(\xi,\xi))(1-
\frac{\Vert\Kxi\Vert^2}{r^2})\\
&&\\
 &&\ \ \ \ \ \ \ \ \ \ \ \ \ \ \ \ \ \ =\frac{r^2}{r^2+\Kxixi^2}g(K(\tilde
e_i,\frac{\xi}{r}),\frac{K(\xi,\xi)}{\Kxixi})(1- \Vert
K(\frac{\xi}{r},\xi)\Vert^2).
\end{array}
\end{equation}

In  both cases: $K(\xi,\xi)=0$ and $K(\xi,\xi)\ne 0$ we have the
following formulas, where each of the indices $i,j$  belongs to
$\{1,...,n\}$ or $\{1,..., n-1\}$, depending on the case. We have
 % Note that $g(\Kxi,e_i)=0$ for $i=1,...,n-1$.

\begin{equation}\label{h_tij}
\begin{array}{rcl}
&&\h_\xi(\tilde e_i^v,e_j^h)=\f\{g(K(\tilde
e_i,e_j),\xi)\\
&&\\
&&\ \ \ \ \ \ \ \ \ \ \ \ \ \ \ \ \ \ \ \ \
-\frac{1}{r^2}g(K(\xi,\xi),e_j)g(K(\xi,\xi),K(\xi, \tilde
e_i))\\
&&\\
&&\ \ \ \ \ \ \ \ \ \ \ \ \ \ \ \ \ \ \ \ \ \ \ \ \ \ \ \ \ \ \ \ \
\ \ \ \
-\frac{1}{2}g(R(e_j,\Kxi)\xi,\tilde e_i)\}\\
&&\\
&&\ \ \ \ \ \ \ \ \ \ \ \ \ \
=\frac{r}{\sqrt{r^2+\Kxixi^2}}\{g(K(\tilde
e_i,e_j),\frac{\xi}{r})\\
&&\\
&&\ \ \ \ \ \ \ \ \ \ \ \ \ \ \ \ \ \ \ \ \ \
-g(K(\frac{\xi}{r},\frac{\xi}{r}),e_j)(K(\xi,\xi),K(\frac{\xi}{r},
\tilde
e_i))\\
&&\ \ \ \ \ \ \ \ \ \ \ \ \ \ \ \ \ \ \ \ \ \ \ \ \ \ \ \ \ \ \ \ \
\ \ \  \ \  -\frac{1}{2}g(R(e_j,\Kxi)\frac{\xi}{r},\tilde e_i)\},
\end{array}
\end{equation}
\bigskip

\begin{equation}\label{h_ij}
\begin{array}{rcl}
&&\h_\xi(e_i^h,e_j^h)=\frac{1}{\sqrt{r^2+\Kxixi^2}}\{\frac{1}{2}[g((\gnabla_{e_i}K)(\xi,\xi),e_j)\\
&&\\
&& \ \ \ \ \ \
+g((\gnabla_{e_j}K)(\xi,\xi),e_i)]-2g(K(\xi,e_i),K(\xi,e_j)
\}\\
&&\\
&&\ \ \ \ \ \ \ \ \ \ =\frac{r}{\sqrt{r^2+\Kxixi^2}}\{\frac{1}{2}[g((\gnabla_{e_i}K)(\frac{\xi}{r},\xi),e_j)\\
&&\\
&& \ \ \ \ \ \ +g((\gnabla_{e_j}K) (\frac{\xi}{r},\xi),e_i)]
-2g(K(\frac{\xi}{r},e_i),K(\xi,e_j)) \},
\end{array}
\end{equation}
\bigskip

\begin{equation}\label{h_titj}
\h_\xi(\tilde e_i^v,\tilde e_j^v)=\f\{-g(\tilde e_i,\tilde
e_j)-g(K(\tilde e_i,\tilde e_j),K(\xi,\xi))\}.
\end{equation}
\bigskip
From the last formula we get
\medskip
\begin{equation}\label{h2_titj}
\begin{array}{rcl}
&&\sum_{i,j=1}^{n-1}(\h_\xi(\tilde e_i^v,\tilde
e_j^v))^2=\frac{n-1}{r^2+\Kxixi^2}\\
&&\\
&&\ \ \ \ \ \ \ \ \ \ \ \ \ \ \ \ \ \ \ \ \ \
+\frac{2(n-1)}{r^2+\Kxixi^2}\sum_{i=1}^{n-1}
g(K(\tilde e_i,\tilde e_i),K(\xi,\xi))\\
&&\\
&&\ \ \ \ \ \ \ \ \ \ \ \ \ \ \ \ \ \ \ \ \ \ \ +
\frac{1}{r^2+\Kxixi^2}\sum_{i,j=1}^{n-1} g(K(\tilde e_i,\tilde
e_j),K(\xi,\xi))^2\\
&&\\
&&\ \ \ \ \ \ \ \ \ \ \ \ \ \ \ \ \ \ \ \ \ \ \ \ =\frac{n-1}{r^2+\Kxixi^2}\\
&&\\
&&\ \ \ \ \ \ \ \ \ \ \ \ \ \ \ \ \ \ \ \ \ \
+\frac{2r^2(n-1)}{r^2+\Kxixi^2}\sum_{i=1}^{n-1}
g(K(\tilde e_i,\tilde e_i),K(\frac{\xi}{r},\frac{\xi}{r}))\\
&&\\
&&\ \ \ \ \ \ \ \ \ \ \ \ \ \ \ \ \ \ \ \ \ \ \ +
\frac{r^2}{r^2+\Kxixi^2}\sum_{i,j=1}^{n-1} (g(K(\tilde e_i,\tilde
e_j),K(\frac{\xi}{r},\xi)))^2.
\end{array}
\end{equation}

The quantity in (\ref{h2_titj}) is the squared norm of $\h$
restricted to the distribution $\mathcal V ^r$, see (\ref{Vr_Hr}).
In order to estimate the function $\h^2$ for small radii we consider
the smooth mapping $(\lambda, \xi)\to \Vert\h_{\lambda \xi}\Vert^2$
defined on the set $(0,1]\times\Sr$. To get the values of $\h$ at
$\lambda \xi\in S^{\lambda r}\M$, it is sufficient to multiply $r$
and $\xi$ by $\lambda$ in the above formulas. Note that the vectors
$\tilde e_1,..., \tilde e_{n-1}$ as well as the vectors
$e_1,...,e_n$ (in the case where $K(\xi,\xi)=0$) or
$e_1,...,e_{n-1}$ (in the case where $K(\xi,\xi)\ne 0$) can be the
same for $\xi$ and $\lambda\xi$.

 We look at the right hand sides of
the formulas (\ref{h_nn})-(\ref{h_ij}), replace $r$ by $\lambda r$,
$\xi$ by $\lambda \xi$  and treat them as functions of two variables
$\lambda$, $\xi$. The right hand side is the last expression in a
formula. At the moment we do not resort to any geometric
interpretation of the functions. One can observe that the image of
each of the function is  bounded in $\R$. For instance, let us first
consider the function in (\ref{h_ij}). Let $P:(0,1]\times \Sr \to $
be the right hand side of (\ref{h_ij}) with $r$  and $\xi$
multiplied by $\lambda$. Take the following well-defined and
continuous function (defined on a compact set)

\begin{equation}
\begin{array}{rcl}
&&P_1:[0,1]\times \Sr\times (S^1\M)^3\ni(\lambda, \xi, u,v,w)\to\\
&&\\
&&\frac{r}{\sqrt{ r^2+\lambda^2\Kxixi^2}}\{-\frac{1}{2}[g((\gnabla_uK)(w,\lambda \xi),v)\\
&&\\
&& \ \ \ \ \ \ +g((\gnabla_{v}K) (w,\lambda\xi),u)]
+2g(K(w,u),K(\lambda\xi,v)) \}\in\R.
\end{array}
\end{equation}
The image of the function $P$ is contained in the image of the
function $P_1$, which is bounded.

 In the same manner we treat the right hand sides in the formulas
 (\ref{h_nn})-(\ref{h_tij}). It cannot be done with (\ref{h_titj}).
But it can be done with the expression for
$\Vert\h_\xi'\Vert^2-\frac{n-1}{r^2+\Kxixi^2}$, where $\h'_\xi$
stands for the restriction of $\h_\xi$ to $\mathcal V^r_\xi$. It
follows from (\ref{h2_titj}) that the image of the function
$$(0,1)\times \Sr\ni (\lambda,\xi)\to
\Vert\h'\Vert^2-\frac{n-1}{\lambda^2r^2+\lambda^4\Kxixi^2}$$ is
contained in a bounded set of $\R$. Since $\{\Kxixi^2, \ \xi\in
\Sr\}$ is bounded, $\frac{n-1}{\lambda^2r^2+\lambda^4\Kxixi^2}$ is
bigger (on the whole $\Sr$) than any given number for  a
sufficiently small $\lambda$.\koniec
\medskip

Using (\ref{mean_curvature}) and the proof of Theorem
\ref{thm_normh} one gets
\begin{thm}\label{H2-h2}
Let $(\M, g,\nabla)$  be a compact statistical manifold whose
dimension is greater than $\rm 2$ and $\h$, $H$ stand for the second
fundamental form and the mean curvature of $\Sr$ in $T\M$. For an
arbitrary number $L$ there exists a sufficiently small radius $r$
such that the function $H^2-\Vert \h\Vert^2$  is grater than $L$
everywhere on $\Sr$.
\end{thm}
\proof Using  (\ref{mean_curvature}) and the same method as in the
proof of Theorem \ref{thm_normh} we see that the
\begin{eqnarray*}
&&(0,1]\times \Sr\ni(\lambda, \xi)\to
\left(H^2_{\lambda\xi}-\frac{(n-1)^2}{\lambda^2r^2+\lambda
^4\Kxixi^2}\right)\\
&&\\
&&\ \ \ \ \ \ \ \ \ \ \ \ \ \ \ \ \ \ \  -\left(\Vert
\h_{\lambda\xi}\Vert
^2-\frac{n-1}{\lambda^2r^2+\lambda^4\Kxixi^2}\right)
\end{eqnarray*}
is bounded. Hence
\begin{equation}
H^2_{\lambda\xi}-\Vert\h_{\lambda\xi}\Vert^2=\frac{(n-1)(n-2)}{\lambda^2r^2+\lambda
^4\Kxixi^2}+P_2(\lambda,\xi),
\end{equation}
where  $P_2$ is some bounded function on $(0,1]\times \Sr$.\koniec

\bigskip

 The following theorem was proved in
\cite{Kow_1}
\begin{thm}
Let $(\M, g)$  be a compact Riemannian manifold and $\dim \M>2$. For
each sufficiently small radius $r>0$ the sphere bundle $\Sr$  has
positive scalar curvature.
\end{thm}

We propose the following generalization of this theorem to the case
of statistical structures
\begin{thm}\label{thm_scalar}
Let $(\M, g,\nabla)$  be a compact statistical manifold and $\dim
\M>2$. For each real number $L$  and sufficiently small radius $r>0$
the sphere bundle $\Sr$ has  scalar curvature greater than $L$.
\end{thm}

\proof The theorem can be proved in the same manner as Theorems
\ref{thm_H}, \ref{thm_normh} and \ref{H2-h2}. By the Gauss equation
for Riemannian hypersurfaces we have

\begin{equation}\label{tilde_rho}
\tilde\rho =\rhotg-2\tRic(\N,\N)+ H^2-\Vert\h\Vert^2,
\end{equation}
where $\tilde\rho$ is the scalar curvature of $\Sr$.

We have explicit  formulas  for the the scalar curvature $\rhotg$,
see (\ref{scalar_for_tg}), for the normal vector field $\N$, see
(\ref{N}), and for the Ricci tensor $\tRic$, see Proposition
\ref{Ricci_tg}. Using the formulas for the first and  second
components of (\ref{tilde_rho}) and replacing $\xi$ by $\lambda\xi$
and $r$ by $\lambda r$ (where $r$ is fixed) and $\lambda\in (0,1]$
one easily sees that the suitable functions of $(\lambda, \xi)\in
(0,1]\times \Sr$ are globally bounded. For instance, let us look at
the second component of (\ref{tilde_rho}). Recall that  $
\N_\xi=\frac{1}{\sqrt{r^2+\Vert K(\xi,\xi)\Vert^2}}(\xi^v_\xi
-K(\xi,\xi)^h_\xi) $. Thus

\begin{equation}\label{tRicNN}
\begin{array}{rcl}
&&-2\tRic_\xi(\N_\xi,\N_\xi) \\
&&\ \ \ \ \ \ =-\frac{2}{r^2+\Vert K(\xi,\xi)\Vert^2}
\tRic_\xi(\xi^v_\xi,\xi^v_\xi)\\
&&\ \ \ \ \ \ \ \ \ \ \ \ \ \ \ \ \ \ +\frac{4}{r^2+\Vert
K(\xi,\xi)\Vert^2}\{\tRic_\xi(\xi^v_\xi,K(\xi,\xi)^h_\xi)\}\\
&&\ \ \ \ \ \  \ \ \  \ \ \ \ \ \ \ \ \ \ \ \ \ \ \ \
 \ \ \ \ \ -\frac{2}{r^2+\Vert K(\xi,\xi)\Vert^2}\tRic_\xi
(K(\xi,\xi)^h_\xi, K(\xi,\xi)^h_\xi).
\end{array}
\end{equation}

We can now observe each component of the above sum separately. By
(\ref{tRic_vv}) we have

\begin{eqnarray*}
 && -\frac{2}{r^2+\Vert K(\xi,\xi)\Vert^2}
\tRic_\xi(\xi^v_\xi,\xi^v_\xi)\\
&&\ \ \   =-\frac{2}{r^2+\Vert K(\xi,\xi)\Vert^2}
\{-2\Ric^K(\xi,\xi)+(\div^\nabla K)(\xi,\xi)-\frac{1}{4}\Vert
R_3(\xi,\xi)\Vert^2\}
\end{eqnarray*}

By replacing $r$ by $\lambda r$ and  $\xi$ by $\lambda\xi$ we get

\begin{eqnarray*}
&& -\frac{2}{(\lambda r)^2+\Vert K(\lambda\xi,\lambda\xi)\Vert^2}
\tRic((\lambda\xi)^v_{(\lambda\xi)},(\lambda\xi)^v_{(\lambda\xi)})\\
&&\ \ \ \ \ \ \ \ =-\frac{2}{\lambda^2(r^2+\lambda^2\Vert
K(\xi,\xi)\Vert^2)} \{-2\Ric^K(\lambda\xi,\lambda\xi)+(\div^\nabla
K)(\lambda\xi,\lambda\xi)\}\\
&&\ \ \ \ \ \ \ \ \ \ \ \ \ \ \ \ \ \ \ \ \ \
+\frac{2}{\lambda^2(r^2+\lambda^2\Vert K(\xi,\xi)\Vert^2)}
\frac{1}{4}\Vert R_3(\lambda\xi,\lambda\xi)\Vert^2
\end{eqnarray*}

\begin{eqnarray*}
&&=-\frac{2}{(r^2+\lambda^2\Vert K(\xi,\xi)\Vert^2)}
\{-2\Ric^K(\xi,\xi)+(\div^\nabla
K)(\xi,\xi)\}\\
&&\ \ \ \ \ \ \ \ \ \ \ \ \ \ \ \ \ \ \ \ \ \
+\frac{\lambda^2}{2(r^2+\lambda^2\Vert K(\xi,\xi)\Vert^2)}\Vert
R_3(\xi,\xi)\Vert^2.
\end{eqnarray*}
The right hand side of this expression is a well-defined continues
function on $[0,1]\times\Sr$, hence it is bounded. One comes to the
same conclusion by observing the two other  components of
(\ref{tRicNN}). It is now sufficient to use Theorem \ref{H2-h2} to
complete the proof. \koniec


\begin{thebibliography}{20}





%\bibitem{CHBY1} Chen B.Y., \emph{Riemannian geometry of Lagrangian
%submanifolds}, Taiwanese . J. Math., 5 (2001) 681-723

\bibitem{D}  Dombrowski P., \emph{On the Geometry of the Tangent Bundle}, J. Reine Angew. Math.
210, (1962), 73-88.






\bibitem{Gud}  Gudmundsson S.,  Kappos  E.,\emph{On the Geometry of Tangent Bundles},
Expo. Math. 20 (2002) 1-41

\bibitem{Kow} Kowalski O., \emph{Curvature of the Induced Riemannian Metric on the
 Tangent Bundle of a Riemannian Manifold}, J. Reine Angew. Math.,
 250 (1971), 124-129

\bibitem{Kow_1} Kowalski O., Sakizawa M., \emph{On tangent sphere
bundles with small or large constant radius}, Ann. Glob. Anal.
Geom., 18, (2000), 207-219




\bibitem{L}  Lauritzen, S.L., \emph{Statistical manifolds}, Siff.
Geom. Stat. Inference, 10, (1987), 163-216



%\bibitem{LSZ}  A-M. Li, U. Simon, G. Zhao, \emph{Global Affine Differential Geometry of
%Hypersurfaces}, Walter de Gruyter,  1993.


 \bibitem{Mat} H. Matsuzoe, J. Inoguchi \emph{Statistical structures on tangent
 bundles}, Appl. Sci. 5, (2003), 55-75


%\bibitem{NS} K. Nomizu, T. Sasaki, \emph{Affine Differential
%Geometry}, Cambridge University Press, 1994.


\bibitem{BW4} Opozda B., \emph{Bochner's technique for statistical
structures}, Ann. Glob. Anal. Geom., 48 (2015), 357-395

\bibitem{seccur} Opozda B., \emph{ A sectional curvature for statistical
structures},  Linear Algebra Appl.,497 (2016), 134-161



\bibitem{Sasaki}  Sasaki S., \emph{ On the differential geometry of tangent bundles of
Riemannian manifolds}, T$\hat{o}$hoku Math. J. 10 (1958), 338-354

\bibitem{Sat}  Satoh H., Almost Hermitian structures on tangent
bundles, (2019), arXiv 1908.10824



\bibitem{shima}  Shima H. \emph{The geometry Hessian structures},
World Scientific Publishing, (2007)





\end{thebibliography}
\end{document}